\documentclass[preprint,11pt]{elsarticle}

\newtheorem{thm}{Theorem}[section]
\newtheorem{lem}{Lemma}[section]

\newtheorem{prop}{Proposition}[section]

\newtheorem{Defi}{Definition}[section]

\newtheorem{rem}[thm]{\bf Remark}

\newcommand{\ignore}[1]{}{}

\usepackage{amssymb}
\usepackage{amsmath}
\usepackage{color}
\usepackage{soul}
\usepackage[dvipsnames]{xcolor}

\usepackage[colorlinks=true, linkcolor=blue, citecolor=blue]{hyperref}

\def\1{{{\mbox{${\rm{1\negthinspace\negthinspace I}}$}}}}

\newcommand{\eref}[1]{(\ref{#1})}

\newcommand\beq{\begin{equation}}
\newcommand\eeq{\end{equation}}

\usepackage{ifthen}
\usepackage{xkeyval}
\usepackage{todonotes}
\begin{document}

\begin{frontmatter}

\title{Deviation inequalities for stochastic approximation by averaging}
\author{Xiequan Fan, \ \   Pierre Alquier \ \   and\ \ Paul Doukhan   }
% \cortext[cor1]{\noindent Corresponding author. \\
%\mbox{\ \ \ \ }\textit{E-mail}: fanxiequan@hotmail.com (X. Fan)\\
%\mbox{\ \ \ \ \ \ \ \ \ \ \  \ \ \ \ } . }
\address{Center for Applied Mathematics, Tianjin University, Tianjin 300072, China.\\
E-mail: fanxiequan@hotmail.com.}
\address{RIKEN AIP, Japan.\\
E-mail: pierrealain.alquier@riken.jp. }
\address{CY University, Cergy-Pontoise, France.\\
E-mail: doukhan@cyu.fr.}

\begin{abstract}
We introduce a class of Markov chains that includes models of stochastic approximation by averaging and non-averaging.
Using a martingale approximation  method, we establish  various deviation inequalities for separately Lipschitz functions of such a chain, with different moment conditions
on some dominating random variables of martingale differences.
Finally, we apply these inequalities to stochastic approximation by averaging  and empirical risk minimisation.
\end{abstract}

\begin{keyword} Deviation inequalities; martingales; iterated random functions; stochastic approximation by averaging;  empirical risk minimisation
\vspace{0.3cm}
\MSC primary 60G42; 60J05;  60F10;   secondary  60E15
\end{keyword}

\end{frontmatter}

% \linenumbers


%\textcolor{magenta}{CORRECTIONS:}
%1. See 2 questions on page 6.

\section{Introduction}
 Markov chains, or iterated random functions, are of fundamental importance to model dependent phenomena. A nice reference on this topic is \cite{DF99}. Probability inequalities for dependent variables were developed in \cite{F73},  and more recently in \cite{LV01,RPR10}  as well as in  \cite{R09,R13,CDM17,D99,FGL12,FGL17}. Most of these papers involve such inequalities for Markov chains. Recently, \cite{DF15} provided such inequalities for contractive Markov chains thanks to a martingale based technique.

% Markov chains or iterated random functions are of a fundamental importance to model dependent phenomena and a nice reference for this is \cite{DF99}. A large amount of probability inequalities  under dependence may be found in the literature, see \cite{F73}  and more recently \cite{LV01}, \cite{RPR10}  as well as in  \cite{R09}, \cite{R13}, \cite{CDM17}, \cite{D99}, \cite{FGL12},   or  \cite{FGL17}. Many papers involve inequalities for Markov chains  and recent martingale based techniques provide reasonable ones for contractive Markov chains as in \cite{DF15}; such contractive Markov chains are weakly dependent.

 In  these papers, only time homogeneous contractive Markov chains are considered. However, in many practical situations, such as stochastic approximation algorithms \cite{PJ92} and unit roots \cite{P87}, the contraction coefficients are time-varying, and will tend either to 0 or to 1 as $n\to \infty$. In this paper, our objective is to provide results for such non-homogeneous Markov chains. Our framework is a large class of non-homogeneous models introduced in Section \ref{secmodel}. Practical examples of chains fitting such conditions are considered  in Section \ref{secexample}.

%The above references mainly correspond to time homogeneous contractive cases, and we aim at proving results for time non-homogeneous Markov chains. This is the setting of the large class of models introduced in Section \ref{secmodel}. Different situations of stochastic algorithms \cite{PJ92} and unit roots \cite{P87} correspond to such varying contraction coefficients tending either to 0 or to 1 as $n\to \infty$. Several relevant models fitting such conditions are considered  in Section \ref{secexample}.

Using the martingale approximation method developed in \cite{DF15}, we establish various deviation inequalities for separately Lipschitz functions of such chains in Section \ref{sec3}. Our inequalities hold under various moment conditions on some dominating random variables of the martingale differences. Section \ref{deviationiq} is dedicated to various classes of $L^p$-norm concentration inequalities, such as Bernstein type inequalities, semi-exponential bound, Fuk-Nagaev inequalities, as well as von Bahr-Esseen, McDiarmid and Hoeffding type bounds. Section \ref{secM} is devoted to moment inequalities: Marcinkiewicz-Zygmund and von Bahr-Esseen type bounds. Finally, in Section \ref{appl} we apply these inequalities to the stochastic approximation by averaging in Subsection \ref{SA} and to empirical risk minimisation (ERM) in Subsection \ref{ERM}.

%Using the martingale approximation method already used in \cite{DF15}, we establish various deviation inequalities for separately Lipschitz functions of such a chain considered in Section \ref{sec3}, with different moment conditions on some dominating random variables of martingale differences; the necessary machinery is also developed. Section \ref{deviationiq} dedicates to various classes of $L^p$-norm concentration inequalities, such as Bernstein type inequalities, semi-exponential bound, Fuk-Nagaev inequalities, as well as von Bahr-Esseen, McDiarmid and Hoeffding type bounds. Section \ref{secM}  devotes to  moment inequalities of Marcinkiewicz-Zygmund  and von Bahr-Esseen type  bounds. Finally, in Section \ref{appl} we apply these inequalities to the stochastic approximation by averaging in Subsection \ref{SA} and to Empirical Risk Minimisation in Subsection \ref{ERM}.

\subsection{Notations}\label{secnotations}
In  the paper, we adopt the convention that each $x\in\mathbb{R}^d$ is a column vector. The entries of $x$ will be denoted by $x^{(1)},\dots,x^{(d)}$. The transpose of $x$ will be denoted by $x^T$, thus $x^T=(x^{(1)},\dots,x^{(d)})$. The set of $d_1\times d_2$ real-valued matrices will be denoted by $\mathbb{R}^{d_1\times d_2}$, and $I_d$ will denote the identity matrix in $\mathbb{R}^{d\times d}$. Let $\|\cdot\|$ denote a norm on $\mathbb{R}^d$. In most cases, we will use $L^p$ norms. In this case, we will explicitely state that $\|\cdot\|=\|\cdot\|_p$, where
 $$\| x \|_\infty =\max_{1\leq i\leq d} |x^{(i)}| \ \ \ \   \textrm{and} \ \ \ \ \| x \|_p=\Big(\sum_{i=1}^d\displaystyle |x^{(i)}|^p \Big)^{1/p},\  \ \  p \in [1, \infty). $$
For any $M\in\mathbb{R}^{d\times d}$, we put
\[
\lambda_{\min}^{(p)}(M) = \inf_{v\neq 0} \frac{\|Mx\|_p}{\|x\|_p} \text{ \ \ \ and \ \ \  }  \lambda_{\max}^{(p)}(M) = \textcolor{red} \displaystyle \sup_{x\neq 0} \frac{\|Mx\|_p}{\|x\|_p} .
\]

Let $(\Omega, {\mathcal A}, {\mathbb P})$ be a probability space. All the random variables in the paper are defined over $(\Omega, {\mathcal A}, {\mathbb P})$. When $V$ is a nonnegative real-valued random variable, we will let $\|V\|_\infty$ denote its essential supremum (note that there will be no ambiguity with the above). Finally, $({\mathcal X}, d)$ and  $({\mathcal Y}, \delta)$ are two complete separable metric spaces. Our non-homogeneous Markov chains will take values in ${\mathcal X}$.

\subsection{A class of iterated random functions}\label{secmodel}
Let $(\varepsilon_i)_{i \geq 1}$ be a sequence of independent copies of   a  ${\mathcal Y}$-valued
random variable $\varepsilon$. Let $X_1$ be a ${\mathcal X}$-valued random variable independent of $(\varepsilon_i)_{i \geq 2}$. We consider the Markov chain $(X_i)_{i \geq 1}$
such that
\begin{equation}\label{Mchain}
X_n=F_n(X_{n-1}, \varepsilon_n), \quad \text{for any $n\geq 2$},
\end{equation}
where $F_n: {\mathcal X} \times {\mathcal Y} \rightarrow {\mathcal X}$
satisfies that there exists a positive number  $n_0 $ such that
for any  $n \geq n_0,$
\begin{equation}\label{contract}
{\mathbb E}\big [ d\big(F_n(x, \varepsilon_1), F_n(x', \varepsilon_1)\big) \big] \leq
\rho_n d(x, x')
\end{equation}
for some $\rho_n \in [0,1)$, and
\begin{equation}\label{c2}
d(F_n(x,y), F_n(x,y')) \leq \tau_n \delta(y,y')+ \textcolor{blue}{\xi_n}
\end{equation}
for some $\tau_n \geq 0$, $\xi_n \geq 0$. The case $\xi_n\equiv 0$ corresponds to functions $F_n$ that are Lipschitz with respect to $\varepsilon_n$, while the case $\tau_n \equiv 0$ corresponds to bounded chains. Note that when $\tau_n\equiv 0$, the metric $\delta$ is not involved in the properties of the chain.

The case where $F_n\equiv F, \rho_n \equiv \rho$, $\tau_n\equiv \tau$  and $ \xi_n\equiv 0$  for two constants $\rho$ and $\tau$ has been studied by Dedecker and Fan \cite{DF15}.
See also Dedecker, Doukhan and  Fan \cite{DDF19} who  weakened  the condition in~\eqref{c2}.
In these papers, the authors have established very precise inequalities for Lipschitz
functionals of the chain, by assuming various moment conditions.
However, the  conditions $F_n\equiv F$ and $\rho_n \equiv \rho$   are restrictive. They are not satisfied  in many extremely useful models.  For instance,  the recursive algorithm of
 stochastic approximation in Polyak and Juditsky \cite{PJ92} returns a chain for which  the conditions \eqref{Mchain}--\eqref{c2} are satisfied with $F_n(x, y)=( I_d - \tau_n A ) x +\tau_n B  - \tau_n y$ and   $\rho_n = 1-c \, \tau_n$ with $ \tau_n \rightarrow 0$,
  where $A \in \mathbb{R}^{d\times d}$ is a positive-definite matrix, $ x, y, B\in\mathbb{R}^{d}$  and $c, \tau_n >0.$
    A  special case of interest corresponds to
   $\rho_n=1-c_1/n^\alpha$ and $\tau_n = c_2/n^\alpha$  for   three positive constants $c_1,$ $c_2$ and $ \alpha \in (0, 1)$. A second class of frequently used models which
 do not satisfy the condition $\rho_n \equiv \rho$ is that of  time series auto-regressions with a \textcolor{blue}{unit  root,}
 see Phillips and Magdalinos \cite{PM07}. In the model of Phillips and Magdalinos \cite{PM07}, the conditions \eqref{Mchain}--\eqref{c2}  are satisfied with
   $\rho_n=1- c_1/n^\alpha$ and $\tau_n = c_2$. See also  Phillips \cite{P87,P88} for the case $\rho_n =e^{-c_1/n}$ and $\tau_n = c_2$.

\subsection{Examples}\label{secexample}
In this subsection, we give a non exhaustive list of models satisfying the conditions \eqref{Mchain}--\eqref{c2}.

\textbf{Example 1.} In the case where ${\mathcal X}$ is a separable Banach space with norm  $\|\cdot\|_{{\mathcal X}}$, and $d(x,x')=\|x-x'\|_{{\mathcal X}}$,  let us consider the following functional auto-regressive model
\begin{eqnarray}\label{standard}
X_n=f(X_{n-1}) + g(\varepsilon_n)\, ,
\end{eqnarray}
where $f : {\mathcal X}  \rightarrow {\mathcal X}$ and $g: {\mathcal Y}  \rightarrow {\mathcal X}$ are such that
$$
  \|f(x)-f(x')\|_{{\mathcal X}} \leq \rho \|x-x'\|_{{\mathcal X}}  \quad \text{and} \quad \|g(y)-g(y')\|_{{\mathcal X}} \leq   \delta(y,y')
$$
for some constant $\rho \in [0, 1).$
In this model, the conditions \eqref{Mchain}--\eqref{c2} are
satisfied with
\begin{eqnarray}\label{standard2}
F_n(x, y)=f( x ) + g(y ),  \ \  \rho_n=\rho, \ \  \tau_n=1  \ \ \textrm{and} \ \  \xi_n \equiv 0
\end{eqnarray}
for any $n \geq 1$. This model is a typical example considered in Dedecker and Fan \cite{DF15}.
We refer to the papers by Diaconis and Freedman \cite{DF99} and Alquier et al. \cite{ABDG20}  for many other interesting examples.

\textbf{Example 2.}  Consider the following auto-regression with a unit root model (see Phillips \cite{P87,P88} or Phillips and Magdalinos \cite{PM07}): for any $n\geq 2$,
\begin{eqnarray}
X_n &=& \frac{1}{1+c/n^\alpha } X_{n-1} +  \varepsilon_n  \quad    \ \ \ \ \textrm{or} \\
X_n &=& (1- \frac{c}{n^\alpha} ) X_{n-1} +  \varepsilon_n,
\end{eqnarray}
where $X_n, \varepsilon_n \in \mathbb{R}, \alpha \in (0, 1)$ and $c$ is a positive constant.
Let $d(x, x')= \delta(x,x')=|x-x'|$.
In this model, the conditions \eqref{Mchain}--\eqref{c2} are
satisfied with $$   \rho_n = 1-   \frac{ c }{n^\alpha}, \ \ \  \tau_n=1 \ \ \  \textrm{and}  \ \ \   \xi_n \equiv 0$$
for any $n$ large enough. Moreover, if $c \in (0, 1),$
the conditions \eqref{Mchain}--\eqref{c2} are
satisfied for any $n \geq 2.$

%
%
%\textbf{Example 3}: Consider the linear problem, see Polyak and Juditsky \cite{PJ92}. We want to find $x^*$, which is  solution of the following equation:
%\begin{equation}
%Ax=B,
%\end{equation}
%where $A \in \mathbb{R}^{d\times d}, x, B\in\mathbb{R}^{d}. $  Assume that $A$ is  positive-definite.
%To obtain the sequence of estimates $(\overline{x}_n)_{n\geq1}$ of the solution $x^*$, the following recursive algorithm will be applied: for $n\geq 2$,
%\begin{eqnarray}
%x_n &=&   x_{n-1} - \frac{\gamma}{n^\alpha} y_n , \quad \quad \quad \quad   y_n\ =\ Ax_{n-1}-B+\varepsilon_n, \\
%\overline{x}_n &=& \frac1n \sum_{i=1}^{n}x_i,
%\end{eqnarray}
%where $\gamma \in (0, \infty),  \alpha \in (0, 1]$ and $x_1$ can be  an arbitrary point or a random  point  in $\mathbb{R}^d.$ The sequence $(y_n)_{n\geq1}$ is observed,
%with $\varepsilon_n\in\mathbb{R}^{d}. $ Here $Ax_{n-1}-B$ is a prediction residual and $\varepsilon_n$ is a random disturbance.
%Let $d$ and $\delta$ be the Euclid norm on $\mathbb{R}^{d}$.
% In this recursive  algorithm, the conditions \eqref{Mchain} - \eqref{c2} are
%satisfied with $$F_n(x, y)=( I -  \frac{\gamma}{n^\alpha}A) x +\frac{\gamma}{n^\alpha}B  - \frac{1}{n^\alpha}y , \ \ \ \  \rho_n = 1- \frac{\gamma \lambda_{\min}(A) }{n^\alpha}   \ \ \  \textrm{and}  \ \ \  \tau_n=\frac{\gamma}{n^\alpha} $$
%for $n$ large enough, where $ \lambda_{\min}(A)$ is the smallest eigenvalues of the matrix $A.$ Moreover, if $\gamma\lambda_{\min}(A) \in (0, 1),$
%the conditions \eqref{Mchain} - \eqref{c2} are
%satisfied for all $n  \geq 2.$

\textbf{Example 3.} Consider the following generalized linear problem.  Set $\mathcal{X}={\mathcal Y}=\mathbb{R}^d$, and let $d$ and $\delta$ be the $L^p$-norm on $\mathbb{R}^{d}$, that is $d(x, x') = \delta(x, x')  =\|x-x'\|_p$ for $p \in [1, \infty]$.  Assume that $(A_i)_{i\geq 1}$  is a  sequence of positive-definite i.i.d.\ random matrices such that $\mathbb{E}A_i=A\in \mathbb{R}^{d\times d}$, $ \lambda^{(p)}_{\min}(A_1) \geq \lambda $ almost surely for some positive constant $\lambda$ and  $\|\lambda^{(p)}_{\max}(A_1) \|_{\infty} < \infty$,  and
that $(B_i)_{i\geq 1}$  is a sequence  of i.i.d. random vectors such that   $\mathbb{E}B_i=B\in\mathbb{R}^{d}$.
 Here, for any given $i$,  $A_i$ and $B_i$ may not be independent. These sequences are observed.
 We want to find $x^*$, which is  solution of the following equation:
\begin{equation}\label{fsdfs12}
Ax=B.
\end{equation}
To obtain the sequence of estimates $(\overline{X}_n)_{n\geq1}$ of the solution $x^*$, the following recursive algorithm will be applied: for any $n\geq 2$,
\begin{eqnarray}
X_n &=&   X_{n-1} - \frac{\gamma}{n^\alpha} Y_n , \quad \quad \quad \quad   Y_n\ =\ A_{n-1}X_{n-1}-B_{n-1}+\eta_n, \\
\overline{X}_n &=& \frac1n \sum_{i=1}^{n}X_i,
\end{eqnarray}
where  $\gamma \in (0, \infty) $  and $\alpha \in [0, 1)$  are constants, and $X_1 \in \mathbb{R}^d$ can be an arbitrary deterministic  point  or a random  point independent  of  $(A_n)_{n\geq1}$ and $(B_n)_{n\geq1}$.    Here $A_{n-1}X_{n-1}-B_{n-1}$ is the prediction residual, and $\eta_n\in\mathbb{R}^{d}$ is a random  disturbance  independent of  $(A_n)_{n\geq1}$ and $(B_n)_{n\geq1}$  (the distinction between $B_i$ and $\eta_i$ is kept because in applications, the user might add a random perturbation $\eta_i$ to  the  noise of the gradient $B-B_i$).

In the special case $A_i\equiv A$ and $B_i\equiv B$ for any $i$, the generalized linear problem  becomes the usual linear problem, see Polyak and Juditsky \cite{PJ92}. More generally, when $B_i$ can be random and $A_i \equiv A$, this example matches our framework. Indeed, put $\varepsilon_i=\eta_i-B_{i-1}$, the conditions \eqref{Mchain}--\eqref{c2} are satisfied with
$$F_n(x,y) = F_n(x,y) =( I_d -  \frac{\gamma}{n^\alpha}A) x  - \frac{\gamma}{n^\alpha}y , \ \ \ \  \rho_n = 1- \frac{\gamma \lambda_{\min}^{(p)}(A)  }{n^\alpha} ,\ \  \tau_n=\frac{\gamma}{n^\alpha} \ \  \textrm{and}    \ \ \xi_n \equiv 0 $$
for any $n \geq 2$.

 In the general case, $A_i$ and $B_i$ are random, so we will enrich the variable $\varepsilon_i$ by $\varepsilon_i=( A_{i-1}, \eta_i-B_{i-1})$. The conditions might not be satisfied in this case. However, it is quite common to seek for the best approximation of $x^*$ in
the set $\mathcal{C}=\{x: \, \| x\|_2 \leq D \},$  that is $x^* \in \mathcal{C}$ and
\begin{equation}
 \|Ax^*-B\|_2=   \min_{x\in\mathcal{C}}\|Ax-B\|_2.
\end{equation}
In this case, it is natural to add a projection step on $\mathcal{C}$.  We then focus on $p=2$.   Let $\Pi_{\mathcal{C}}: \mathbb{R}^{d} \rightarrow \mathcal{C}$ denote the orthogonal projection on $\mathcal{C}$.   Note that $\Pi_{\mathcal{C}}$ is
such that $\|\Pi_{\mathcal{C}} x - \Pi_{\mathcal{C}} y \|_2 \leq \|x -y\|_2$ for any $ x, y \in \mathbb{R}^d$ and  $\Pi_{\mathcal{C}} x =x$ for any $x \in \mathcal{C}$. Then, for any $n\geq 2$, take
\begin{eqnarray}
X_n &=& \Pi_\mathcal{C} \Big[X_{n-1} - \frac{\gamma}{n^\alpha} Y_n \Big], \quad \quad \quad \quad   Y_n\ =\ A_{n-1}X_{n-1}-B_{n-1}+\eta_n, \\
\overline{X}_n &=& \frac1n \sum_{i=1}^{n}X_i.
\end{eqnarray}
The conditions \eqref{Mchain}--\eqref{c2} are satisfied, for $y=(M, u)$, with
$$F_n(x, y)= \Pi_{\mathcal{C}} \left[ ( I_d -  \frac{\gamma}{n^\alpha}M)  \Pi_{\mathcal{C}}[ x ]   - \frac{\gamma}{n^\alpha} u  \right] ,   \ \   \rho_n = 1- \frac{\gamma \lambda   }{n^\alpha},  \ \    \tau_n= \frac{\gamma}{n^{\alpha}} \  \ \textrm{and}\ \  \xi_n= \frac{2 D \gamma\| \lambda_{\max}^{(2)}(A_1)   \|_{\infty} }{n^{\alpha}}  $$
for any $n \geq 3$  and $\alpha \in (0, 1)$. For the case $\alpha =0,$  the conditions \eqref{Mchain}--\eqref{c2} are also satisfied, but with an additional assumption that  $ \gamma  < 1/ \lambda  $.

\textbf{Example 4.} For the usual linear problem (cf.\ equation (\ref{fsdfs12})), another  recursive algorithm may be applied: for any $n\geq 2$,
\begin{eqnarray}
X_n &=&   X_{n-1} -\gamma \, Y_n , \quad \quad \quad \quad   Y_n\ =\ A X_{n-1}-B +\frac{1}{n^\alpha}\varepsilon_n, \\
\overline{X}_n &=& \frac1n \sum_{i=1}^{n}X_i,
\end{eqnarray}
where  $\alpha \in (0, 1)$, $\gamma \in(0, \infty)$ is a constant  such that   $\gamma  \lambda_{\min}^{(p)} (A) \in (0, 1)$, and $X_1 \in \mathbb{R}^d$ can be an arbitrary deterministic point or a random  point independent  of  $(\varepsilon_n)_{n\geq 2}$.
Let $d$ and $\delta$ be the $L^p$-norm   on $\mathbb{R}^{d}$.
 In this recursive  algorithm, the conditions \eqref{Mchain}--\eqref{c2} are
satisfied with $$F_n(x, y)=( I_d -  \gamma   A ) x + \gamma B  - \frac{\gamma}{n^\alpha}y , \ \ \ \  \rho_n = 1-  \gamma \lambda_{\min}^{(p)}(A),     \ \  \tau_n=\frac{\gamma}{n^\alpha}   \ \  \textrm{and}  \ \ \, \xi_n \equiv 0  $$
for any $n \geq 2$.

\textbf{Example 5.}   A  third recursive algorithm for the usual linear problem  is given by: for any $n\geq 2$,
\begin{eqnarray}
X_n &=&   X_{n-1} -\frac{\gamma}{n^\alpha} \, Y_n+\varepsilon_n , \quad \quad \quad \quad   Y_n\ =\ A X_{n-1}-B , \\
\overline{X}_n &=& \frac1n \sum_{i=1}^{n}X_i,
\end{eqnarray}
where  $\alpha \in (0, 1)$, $\gamma \in(0, \infty)$ are constants such that   $\gamma  \lambda^{(p)}_{\min}(A)  \in (0, 1)$, and $X_1 \in \mathbb{R}^d$ can be an arbitrary deterministic point or a random  point independent of $(\varepsilon_n)_{n\geq 2}$.
Let $d$ and $\delta$ be the $L^p$-norm   on $\mathbb{R}^{d}$.
 In this recursive  algorithm, the conditions \eqref{Mchain}--\eqref{c2} are
satisfied with $$F_n(x, y)=( I_d - \frac{\gamma}{n^\alpha}   A ) x + \frac{\gamma}{n^\alpha} B   -  \gamma y , \ \ \ \  \rho_n = 1-  \frac{\gamma \lambda_{\min}^{(p)}(A)  }{n^\alpha},     \ \  \tau_n\equiv \gamma  \ \  \textrm{and}  \ \    \xi_n \equiv 0    $$
for any $n \geq 2$.\bigskip

\textbf{Example 6.} We extend the previous examples to optimization of non-linear functions. We still consider $\mathcal{X}={\mathcal Y}=\mathbb{R}^d$  and   focus on the $L^2$ norm in this example. In machine learning, we need to minimize a function involving a large number of differentiable terms $L(x)=\sum_{i=1}^N \ell_{i}(x) $ on the set $\mathcal{C}=\{x: \, \| x\|_2 \leq D \}$.  A popular strategy to this end  is to use the projected stochastic gradient descent (SGD): for any $n\geq 2,$
$$ X_{n} = \Pi_\mathcal{C} \left[ X_{n-1} - \frac{\gamma}{n^{\alpha}} \hat{\nabla}_{J_n} L(X_{n-1}) \right], $$
where $\alpha \in (0, 1]$,  $\Pi_{\mathcal{C}}: \mathbb{R}^{d} \rightarrow \mathcal{C}$ denote the orthogonal projection on $\mathcal{C}$,  $ J_n $   is drawn uniformly among all the subsets of $\{1, \dots, N\}$ with cardinality $M$ and
$$ \hat{\nabla}_{ J_n } L(x) := \frac{1}{M}\sum_{i\in  J_n } \nabla \ell_{i}(x) .$$
Note that $\mathbb{E}[ \hat{\nabla}_{ J_n } L(x)] = \nabla  L(x) $ for any $x$. More generally, the stochastic gradient Langevin descent (SGLD) is given by
$$ X_{n} = \Pi_\mathcal{C} \left[ X_{n-1} - \frac{\gamma}{n^{\alpha}} \hat{\nabla}_{ J_n } L(X_{n-1}) - \frac{\gamma}{n^{\alpha}}\eta_{n} \right] $$
for some i.i.d.\ sequence $\eta_n$ of random perturbations added by the user.
For $y=(J, u)$, define $$ F_n(x,y) = \Pi_\mathcal{C} \Big[  x - \frac{\gamma}{n^{\alpha}} \hat{\nabla}_{ J }  L(x) -  \frac{\gamma}{n^{\alpha}} u  \Big]. $$
  If we define $\varepsilon_i=(J_i, \eta_i)$, then this example fits \eqref{Mchain}. It  is easy to see that
\begin{align*}
 &\| F_n(x,y) - F_n(x',y) \|_2^2
  \ \leq \   \|  x - \frac{\gamma}{n^{\alpha}} \hat{\nabla}_{J_n}L (x) - x' - \frac{\gamma}{n^{\alpha}} \hat{\nabla}_{J_n}L (x') \|_2^2
  \\
 &\quad\quad\quad\quad  = \|x-x'\|_2^2 + \Big\| \frac{\gamma}{n^{\alpha}} \hat{\nabla}_{J_n} L(x) - \frac{\gamma}{n^{\alpha}} \hat{\nabla}_{J_n}L (x')  \Big\|_2^2
   - 2 (x-x')^T \left(\frac{\gamma}{n^{\alpha}} \hat{\nabla}_{J_n} L (x) - \frac{\gamma}{n^{\alpha}} \hat{\nabla}_{J_n} L(x') \right).
\end{align*}
Common assumptions are that the $\ell_i$'s are $m$-strongly convex, $m > 0,$  which gives
$$ (x-x')^T \left(\hat{\nabla}_{J_n} L(x) -  \hat{\nabla}_{J_n} L(x') \right) \geq m \|x-x'\|_2^2  $$
and that  their gradients  are $\ell$-Lipschitz,  that is,
$$ \| \hat{\nabla}_{ J_n }L (x) -  \hat{\nabla}_{ J_n } L(x') \|_2^2 \leq \ell^2 \|x-x'\|_2^2. $$
  We then obtain
\begin{align*}
\| F_n(x,y) - F_n(x',y) \|_2^2 &\leq   \left(1-\frac{2m\gamma}{n^{\alpha}}  + \frac{\ell^2\gamma^2}{n^{2 \alpha}} \right) \|x-x'\|_2^2
 \\
 &  =  \left[  \left(1-\frac{m\gamma}{n^\alpha}\right)^2 + \frac{(\ell^2-m^2)\gamma^2}{n^{2\alpha}} \right] \|x-x'\|_2^2
\end{align*}
and thus the condition (\ref{contract}) is satisfied with
$$
\rho_n = \sqrt{ \left(1-\frac{m\gamma}{n^\alpha}\right)^2 + \frac{(\ell^2-m^2)\gamma^2}{n^{2\alpha}}}.
$$
Note that $\rho_n- 1\sim m\gamma/n^{\alpha}$. So for any $n$ large enough, we have for example
$$ \rho_n \leq 1 - \frac{m\gamma}{2n^{\alpha}}. $$
The  condition also holds  in the case $\alpha =0,$ but with an additional assumption that  $ 2m\gamma   - \ell^2\gamma^2 \in (0, 1)$, that can always be achieved with an adequate choice of $\gamma$.  Finally, let us assume that
 $\| \nabla \ell_i(x) \|_2 \leq B$ for any $x \in \mathcal{C}$ and some $B>0$.
Condition \eqref{c2} is satisfied with
$$  \tau_n= \frac{\gamma}{n^{\alpha}}  \  \ \textrm{and}  \  \ \xi_n= \frac{2 B \gamma  }{n^{\alpha}}  $$
for any $n$ large enough. 

\textbf{Example 7.}  Our final example illustrates that non-homogeneity can appear even in the context of a time homogeneous chain, if it is only observed at non evenly spaced dates $t_1,t_2,\ldots.$ Assume $(t_i)_{i\geq1}$ is an increasing sequence, put $k_1=t_1$ and $k_i=t_i-t_{i-1}>0$ for any $i>1$. Consider $F_n \equiv F$, $\rho_n \equiv \rho$, $\tau_n\equiv \tau$,  $ \xi_n \equiv 0 $  and $(X_i)_{i\geq 1}$ the corresponding chain in~\eqref{Mchain}--\eref{c2}. Assume that only the subsequence $(X_{t_i})_{i\geq 1}$ is observed. Let $\varepsilon^{(Z)}_i=(\varepsilon_{j,i})_{j\ge 1}$ be an i.i.d. copy of the sequence  $(\varepsilon_i)_{i\geq 1}$, and define, for any $n\geq 1$ and $y=(y_i)_{i\geq 1}$,
$$F^{(Z)}_n(x,y)=F(F(\ldots F(F(x,y_1),y_2),\ldots , y_{k_{n}-1}),y_{k_n}) $$
($F_n$ only depends on the $k_n$ first terms of the sequence $y$).
It is clear that $Z_n=F^{(Z)}_n(Z_{n-1},\varepsilon^{(Z)}_n)$ admits the same distribution as $X_{t_n}$. Then with the notations above, $\rho_n^{(Z)}=\rho^{k_n}$, and this quantity tends to $0$ as $n\to\infty$ if $k_n\to\infty$, which corresponds to a situation where sampling times become rarer and rarer. The expression of $\tau_n^{(Z)}$ is not clear in general. Let us now restrict our attention to the additive model in~\eqref{standard} (but note that even if  \eqref{standard} holds for $(X_i)_{i\geq 1}$, in general a similar expression does not hold for $(Z_i)_{i\geq 1}$). In this case,
\begin{eqnarray}
| F^{(Z)}_n(x,y) - F^{(Z)}_n(x,y') | \leq \tau \sum_{i=1}^n \rho^{n-i} \delta(y_{i},y_{i}'),
\end{eqnarray}
so, for example with the sup metric on ${\cal Y}^{\mathbb{N}_*}$, given by $\sup_{i\in\mathbb{N}_*} \delta(y_i,y_i') $, we obtain $\tau_n^{(Z)} = \frac{\tau(1-\rho^n)}{1-\rho} \leq  \frac{\tau}{1-\rho} $.

\section{Lipschitz functions of random vectors $X_1, \ldots , X_n$ }\label{sec3}
\setcounter{equation}{0}
We remind that  $\| \cdot \|$ is  a norm on $\mathbb{R}^d$. Let $f: {\mathcal X}^n \mapsto {\mathbb R}^d $ be a separately Lipschitz function, that is
\begin{equation} \label{codiMD}
 \| f(x_1, x_2, \ldots, x_n)-f(x'_1, x'_2, \ldots, x'_n)  \|  \leq d(x_1,x'_1)+ d(x_2,x'_2)+\cdots +   d(x_n, x'_n) .
\end{equation}
Let $\mathbb{P}_{X_1}$ and $\mathbb{P}_\varepsilon$  be the distributions of $X_1$ and $\varepsilon$, respectively.
Assume that $\| \cdot \|$   satisfies
\begin{equation} \label{codint}
\Big \|\int h(x)  \mathbb{P}_{X_1} (dx)  \Big \| \leq  \int \| h(x)  \|  \mathbb{P}_{X_1} (dx)   \ \  \   \textrm{and} \ \ \ \Big \|\int h(x)  \mathbb{P}_{\varepsilon} (dx)  \Big \| \leq  \int \| h(x)  \|  \mathbb{P}_{\varepsilon} (dx)
 \end{equation}
 for any measurable function $h: {\mathcal X}^n \mapsto {\mathbb R}^d.$
 Clearly,  if  $\|  \cdot \|= \|  \cdot \|_p, p \in [1, \infty]$,
 then the condition (\ref{codint}) is satisfied.

Let
\begin{equation}\label{Sn}
S_n:=f(X_1, \ldots , X_n) -{\mathbb E}[f(X_1, \ldots , X_n)]\, .
\end{equation}
Denote    $(\mathcal{F}_k)_{k\geq 0}$ the natural filtration of the  chain $(X_k)_{  k \geq1}$, that is  ${\mathcal F}_0=\{\emptyset, \Omega \}$
and for any $k \in {\mathbb{N}}^{*}$,
${\mathcal F}_k= \sigma(X_1, X_2,  \ldots, X_k)$. For any $k \in [0, n]$,
define
\begin{equation}\label{gk}
g_k(X_1, \ldots , X_k)= {\mathbb E}[f(X_1, \ldots, X_n)|{\mathcal F}_k]\,
\end{equation}
and  for any $k \in [1, n]$,
\begin{equation}\label{dk}
M_k=g_k(X_1, \ldots, X_k)-g_{k-1}(X_1, \ldots, X_{k-1})\, .
\end{equation}
Then  $(M_k, \mathcal{F}_k)_{1\leq k \leq n}$ is a finite sequence of martingale differences.
For any $k \in [1, n-1]$, let then
$$
S_k:=M_1+M_2+\cdots + M_k \, ,
$$
  and  note that $S_n$ is already introduced in
\eref{Sn} and satisfies $S_n =M_1+M_2+\cdots + M_n \,$.
Then $(S_k, \mathcal{F}_k)_{1\leq k \leq n}$ is a martingale.

The following proposition gives some interesting properties of
the functions  $(g_k)_{1\leq k \leq n}$ and of the martingale differences $(M_k, \mathcal{F}_k)_{1\leq k \leq n}$.  In this paper, we focus on the case  $n_0=2,$
where $n_0$ is given by the conditions (\ref{contract}) and (\ref{c2}).

\begin{lem}\label{McD}
For any $k \in [1, n]$
and $\rho_k$ in $[0, 1)$, let  $$K_{k,n} =  1 + \rho_{k+1} +  \rho_{k+1}\rho_{k+2}+ \cdots +  \rho_{k+1}\rho_{k+2}\cdots  \rho_{n },\ \ \  k \in [1, n-1] \ \
\textrm{and} \ \   K_{n,n}=1. $$
Let $(X_i)_{i \geq 1}$ be a Markov chain satisfying \eref{Mchain} for some
functions $(F_n)_{n\geq 1}$ satisfying \eref{contract}. We also assume that $\| \cdot \|$   satisfies \eref{codint}. Let $g_k$ and $M_k$ be defined
by
\eref{gk} and \eref{dk}, respectively.
\begin{enumerate}
\item
The function $g_k$ is separately Lipschitz and such that
$$
  \| g_k(x_1, x_2, \ldots, x_k)- g_k(x'_1, x'_2, \ldots, x'_k)  \|  \leq d(x_1, x'_1) +\cdots + d(x_{k-1}, x'_{k-1}) + K_{k,n} d( x_k,  x'_k).
$$
\item  Let $G_{X_1}$ and $H_{k,\varepsilon}$ be  functions defined by
$$
G_{X_1}(x)=\int d(x, x') \mathbb{P}_{X_1}(dx')
$$
and
$$
H_{k,\varepsilon}(x,y)= \int d(F_k(x,y), F_k(x,y'))\mathbb{P}_{\varepsilon}(dy') \, , \ \ k \in [2,n],
$$
respectively. Then, the martingale difference $M_k$ satisfies that
$$
   \|M_1\| \leq K_{ 1,n} G_{X_1}(X_1)
$$
and for any $k \in [2,n]$,
$$
   \|M_k\|\leq  K_{k,n}  H_{k,\varepsilon}(X_{k-1}, \varepsilon_k) \, .
$$
\item
Assume moreover that $F_n$ satisfies \eref{c2}, and
let $G_{\varepsilon}$ be the  function defined by
$$
  G_{\varepsilon}(y)= \int  \delta(y, y')\mathbb{P}_{\varepsilon}(dy')  \, .
$$
Then $H_{k,\varepsilon}(x,y) \leq \tau_k G_\varepsilon(y)$,
and consequently, for any $k \in [2,n]$,
$$
 \|M_k\|\leq   K_{k,n} [ \tau_k  G_\varepsilon(\varepsilon_k)  +  \xi_k  ] \, .
$$

\item
Assume moreover that there exist  three constants $ \alpha \in [0, 1) ,$   $ \rho \in (0, 1)$  and $\eta \in (0, \infty)$  such that
for any $n \geq 2,$
\begin{equation}\label{label1.5}
  \rho_{n}\leq 1- \rho/n^\alpha  \ \ \ \textrm{and }\ \ \  \max\{ \xi_n, \tau_n \} \leq \eta/n^\alpha. \
\end{equation}
Then $ K_{ 1,n}=O(1)$ and   $(K_{ k,n} [\tau_k+\xi_k])_{k\geq 1} $   is uniformly bounded for all $k$ and $n$.

\item
Assume moreover  that there exist three constants $ \alpha \in (0, 1) ,$   $ \rho \in (0, 1)$  and $\eta \in (0, \infty)$  such that
for any  $n \geq 2,$
\begin{equation}\label{label1.6}
   \rho_{n}\leq 1- \rho/n^\alpha  \ \ \ \textrm{and }\ \ \  \max\{ \xi_n, \tau_n \}   \leq \eta.  \ \ \
\end{equation}
Then $ K_{ 1,n}=O(1)$ and  $ K_{ k,n} [\tau_k+\xi_k]= O(k^\alpha)$   as $k \rightarrow \infty.$

\item
Assume moreover that there exist  three  constants $ \alpha \in (0, 1] ,$   $ \rho \in (0, 1)$  and $\eta \in (0, \infty)$  such that
for any  $n \geq 2,$
\begin{equation}\label{label1.7}
   \rho_{n}\leq \rho  \ \ \ \textrm{and }\ \ \  \max\{ \xi_n, \tau_n \} \leq \eta/n^\alpha.   \ \ \
\end{equation}
Then $ K_{ 1,n}=O(1)$ and $ K_{ k,n} \tau_k= O(k^{-\alpha})$ as $k \rightarrow \infty.$

%In particular, it holds
%$$
%   |d_1| \leq  G_{X_1}(X_1)/(1-\rho)\quad \text{and} \quad |d_k|\leq C  %G_\varepsilon(\varepsilon_k)/(1-\rho)\, .
%$$
\end{enumerate}
\end{lem}

\begin{rem}
Let us comment  on the point 4 of Lemma \ref{McD}.
If $\xi_n \equiv \alpha =  0$,  then  $K_{k,n} \leq \sum_{i=0}^{n-k}(1-\rho)^i <\frac{1}{ \rho}$ and $\tau_k\leq \eta$ for any $k \in [1, n]$ and $n$. Thus  $(K_{ k, n} \tau_k)_{k\geq 1} $ is uniformly bounded for all $k$ and $n$, which has been proved by  Proposition 2.1 of Dedecker and Fan \cite{DF15}.
\end{rem}

\begin{rem}
Let us return to the examples in Subsection \ref{secexample}.
It is easy to see that Examples 1  and 7  satisfy the condition (\ref{label1.5}) with $\alpha=0$.  Examples 2 and 5 satisfy the condition (\ref{label1.6}). Examples 3  and 6  satisfy the condition (\ref{label1.5}). Example 4    satisfies the  condition (\ref{label1.7}).
\end{rem}

\noindent {\emph {Proof}.} The first point will be proved by recurrence in the backward sense.
The result is obvious for $k=n$, since $g_n=f$. Assume that it is true at step $k \in [2, n]$, and let us prove it
at step $k-1$. By definition
$$
g_{k-1}(X_1, \ldots, X_{k-1})={\mathbb E}[g_k(X_1, \ldots, X_k)|{\mathcal F}_{k-1}]= \int g_k(X_1, \ldots, X_{k-1}, F_k(X_{k-1},y)) \mathbb{P}_{\varepsilon}(dy)\, .
$$
By assumption (\ref{codint}), it follows that
%\begin{multline}\label{triv1}
%|g_{n-1}(x_1, x_2, \ldots, x_{n-1})-g_{n-1}(x'_1, x'_2, \ldots, x'_{n-1})|\\ \leq
%\int |g_{n }(x_1, x_2, \ldots, F_n(x_{n-1},y))-g_{n }(x'_1, x'_2, \ldots, F_n(x'_{n-1},y))|
%\mathbb{P}_{\varepsilon}(dy)\\
%\leq d(x_1,x'_1)+\cdots + d(x_{n-1}, x'_{n-1})+    \int d(F_n(x_{n-1},y), F_n(x'_{n-1},y))\mathbb{P}_{\varepsilon}(dy)
%\\
%\leq d(x_1,x'_1)+\cdots + d(x_{n-1}, x'_{n-1}) +  \rho_n d(x_{n-1}, x'_{n-1}) \\
%\leq d(x_1,x'_1)+\cdots + d(x_{n-2}, x'_{n-2}) + (1+  \rho_n) d(x_{n-1}, x'_{n-1}) \\
% \leq d(x_1,x'_1)+\cdots + d(x_{n-2}, x'_{n-2}) + K_n  d(x_{n-1}, x'_{n-1})  \, .
%\end{multline}
%\begin{multline}
%|g_{n-2}(x_1, x_2, \ldots, x_{n-2})-g_{n-2}(x'_1, x'_2, \ldots, x'_{n-2})|\\ \leq
%\int |g_{n-1}(x_1, x_2, \ldots, F_{n-1}(x_{n-2},y))-g_{n-1}(x'_1, x'_2, \ldots, F_{n-1}(x'_{n-2},y))|
%\mathbb{P}_{\varepsilon}(dy)\\
%\leq d(x_1,x'_1)+\cdots + d(x_{n-2}, x'_{n-2})+   (1+  \rho_n) \int d(F_{n-1}(x_{n-2},y), F_{n-1}(x'_{n-2},y))\mathbb{P}_{\varepsilon}(dy)
%\\
%\leq d(x_1,x'_1)+\cdots + d(x_{n-3}, x'_{n-3}) + (1+ \rho_{n-1} + \rho_{n-1} \rho_n  ) d(x_{n-2}, x'_{n-2}) \\
%\leq d(x_1,x'_1)+\cdots + d(x_{n-3}, x'_{n-3}) + K_{n-1} d(x_{n-1}, x'_{n-1}) \, .
%\end{multline}
\begin{multline}
  \| g_{k-1}(x_1, x_2, \ldots, x_{k-1})-g_{k-1}(x'_1, x'_2, \ldots, x'_{k-1})  \| \\
  =
\Big\| \int   g_{k}(x_1, x_2, \ldots, x_{k-1}, F_{k}(x_{k-1},y))-g_{k}(x'_1, x'_2, \ldots,x_{k-1}', F_{k}(x'_{k-1},y))
\mathbb{P}_{\varepsilon}(dy) \Big \|\\
\leq
\int   \| g_{k}(x_1, x_2, \ldots, x_{k-1}, F_{k}(x_{k-1},y))-g_{k}(x'_1, x'_2, \ldots,x_{k-1}', F_{k}(x'_{k-1},y))   \|
\mathbb{P}_{\varepsilon}(dy)\\
\leq d(x_1, x'_1) +\cdots + d(x_{k-1}, x'_{k-1}) +   K_{k,n} \int d( F_{k}(x_{k-1},y), F_{k}(x'_{k-1},y) )  \mathbb{P}_{\varepsilon}(dy)
\ \ \ \ \ \\
\leq  d(x_1, x'_1) +\cdots +  d(x_{k-2}, x'_{k-2}) + (1+ \rho_{k} K_{k,n}  ) d(x_{k-1}, x'_{k-1}) \ \ \ \ \ \ \ \ \ \  \ \ \ \ \ \ \ \ \ \ \ \ \ \    \\
\leq  d(x_1, x'_1)+\cdots +  d(x_{k-2}, x'_{k-2}) +  K_{k-1,n} d(x_{k-1}, x'_{k-1}), \ \ \ \ \ \ \ \ \ \ \ \ \ \ \ \ \ \ \ \ \ \ \ \ \ \ \ \ \ \ \ \ \ \ \
\end{multline}
which completes the proof of the point 1.

Let us give a proof for the point  2. First note that
\begin{align}
 \|M_1 \|=\Big \|g_1(X_1)-\int g_1(x)\mathbb{P}_{X_1}(dx)\Big \|\leq
 K_{ 1,n} \int d(X_1, x)\mathbb{P}_{X_1}(dx)=K_{ 1,n}  G_{X_1}(X_1) \, . \label{dfsaf01}
\end{align}
In the same way, for any $k \in [2, n]$,
\begin{align}
  \| M_k  \|  & = \big  \| g_k(X_1, \cdots, X_k)-{\mathbb E}[g_k(X_1, \cdots, X_k)|{\mathcal F}_{k-1}]\big  \| \nonumber \\
 &\leq  \int \big   \| g_k(X_1, \cdots, F_k(X_{k-1}, \varepsilon_k))-g_k(X_1, \cdots, F_k(X_{k-1}, y))\big  \|  \mathbb{P}_{\varepsilon}(dy) \nonumber \\
 &\leq K_{k,n} \int  d(F_k(X_{k-1},\varepsilon_k),
 F_k(X_{k-1}, y)) \mathbb{P}_{\varepsilon}(dy)= K_{ k,n}
  H_{k,\varepsilon}(X_{k-1},\varepsilon_k) \, .  \label{dfsaf02}
\end{align}

 The point  3 is clear, since if \eref{c2} is true, then
$$
H_{k,\varepsilon}(x,y)= \int d( F_k(x,y), F_k(x,y')) \mathbb{P}_{\varepsilon}(dy')
\leq  \int  \tau_k \delta(y, y') \mathbb{P}_{\varepsilon}(dy')   +\xi_k  = \tau_k G_\varepsilon(y)  +\xi_k. 
$$
Consequently, for any $k \in [2,n]$,
$$
  \| M_k \|\leq   K_{ k,n} [ \tau_k  G_\varepsilon(\varepsilon_k) +\xi_k ]  . 
$$

Next we give a proof for the point  4.
%Assume that $\rho_{n}= 1- \eta/n^\alpha $ and $\tau_n=\eta/n^\alpha$ for two constants $\eta>0$ and $ \alpha \in (0, 1).$
It is easy to see that
\begin{eqnarray}
\ln ( \rho_{k+1}\cdots  \rho_{k+l})&=&   \ln \rho_{k+1} + \cdots + \ln  \rho_{k+l} \nonumber \\
   &\leq& - \rho  \sum_{i=k+1}^{k+l} i^{-\alpha}= - \rho \, n^{1-\alpha} \sum_{i=k+1}^{k+l} (\frac i n)^{-\alpha}  n^{-1}  \nonumber \\
   &\leq& - \frac{\rho}{ 1-\alpha } n^{1-\alpha} \Big( ( \frac{ k+l} n)^{1-\alpha}  -  ( \frac{ k+1} n)^{1-\alpha} \Big)  \nonumber \\
   &\leq& - \frac{\rho}{ 1-\alpha } n^{1-\alpha} \frac{1}{1-\alpha} \frac {l-1}n  \frac {n^\alpha}{(k+l)^\alpha} \nonumber \\
   &=& - \frac{\rho}{ (1-\alpha)^2 } \frac {l-1} {(k+l)^\alpha} .  \label{dfsfds}
\end{eqnarray}
Thus, we deduce that
\begin{eqnarray*}
K_{ k,n} \tau_k \!\!  &\leq &\!\!   \frac{\eta}{k^\alpha} \sum_{l= 1}^{n-k }  \exp\Big\{ - \frac{\rho}{ (1-\alpha)^2 }\frac {l-1} {(k+l)^\alpha}  \Big\} \\
   & \leq& \frac{\eta}{k^\alpha} \bigg(  \sum_{l= 1}^{k }  \exp\Big\{ - \frac{\rho}{ (1-\alpha)^2 }\frac {l-1} {(k+l)^\alpha}  \Big\} +\sum_{l= k}^{\infty }  \exp\Big\{ - \frac{\rho}{ (1-\alpha)^2 }\frac {l-1} {(k+l)^\alpha}  \Big\}  \bigg )   \\
   &=:& \frac{\eta}{k^\alpha}  (I_1 +I_2) .
\end{eqnarray*}
For $I_1,$ we have the following estimation
\begin{eqnarray*}
  I_1 \!\!  &\leq &\!\!  \exp\Big\{  \frac{\rho}{ (1-\alpha)^2 }\frac {  1} {(k+1)^\alpha}  \Big\} \sum_{l= 1}^{k }  \exp\Big\{ - \frac{\rho}{ 2^\alpha (1-\alpha)^2 }\frac {l } {k } k^{1-\alpha } \Big\} \\
  &\leq &\!\! k \exp\Big\{  \frac{\rho}{ (1-\alpha)^2 }\frac {  1} {(k+1)^\alpha}  \Big\} \int_{0}^{1 }  \exp\Big\{ - \frac{\rho k^{1-\alpha }}{ 2^\alpha (1-\alpha)^2 } x \Big\}dx \\
  &= &\!\! k \exp\Big\{  \frac{\rho}{ (1-\alpha)^2 }\frac {  1} {(k+1)^\alpha}  \Big\}   \frac{ 2^\alpha (1-\alpha)^2 }{\rho k^{1-\alpha }} \Big (1-\exp\Big\{ - \frac{\rho k^{1-\alpha }}{ 2^\alpha (1-\alpha)^2 }  \Big\} \Big ) \\
   &\leq &\!\! k^\alpha \exp\Big\{  \frac{\rho}{ (1-\alpha)^2 }  \Big\}   \frac{ 2^\alpha (1-\alpha)^2 }{\rho  }  .
\end{eqnarray*}
It is obvious that $I_2$ is bounded for $\alpha \in [0, 1)$.  Hence
$ (K_{ k,n} \tau_k)_{k\geq 1} $ is uniformly bounded for all $k$ and $n$ and $ K_{ 1,n}=O(1)$.

For  the point  5, from the proof of  the point  4, it is easy to see that $ K_{ 1,n}=O(1)$ and  $K_{ k,n} \tau_k =O( k^\alpha), k\rightarrow \infty.$

For the  point  6, it is easy to see that  $K_{k,n} \leq \sum_{i=0}^{n-k}\rho^i <\frac{1}{1-\rho}$ and $\tau_k\leq \eta/k^\alpha$ for all $ k \in [1, n]$. Thus $ K_{ 1,n}=O(1)$ and $  K_{ k, n} \tau_k  =O(1/k^\alpha)$ as $k \rightarrow \infty$. The proof of the lemma is now complete.
\hfill\qed

\section{Deviation inequalities for  $S_n$ with $L^p$-norm} \label{deviationiq}
\setcounter{equation}{0}

Let $n\geq 2.$
In this section, we are interested in the concentration properties of $S_n$ under the $L^p$-norm $\| \cdot \|_p$,
where  $(X_i)_{i \geq 1}$ is a Markov chain satisfying \eref{Mchain}  for some
functions $(F_n)_{n\geq1}$ satisfying \eref{contract} and \eref{c2}.
Clearly, it holds for any $p \in [1, \infty],$
\begin{eqnarray}\label{fsn02}
  \| x \|_p =\Big(\sum_{i=1}^d  | x^{(i)}  |^p \Big)^{1/p} \leq \Big(\sum_{i=1}^d \| x \|_\infty^p \Big)^{1/p}=d^{1/p}\| x \|_\infty. 
\end{eqnarray}
When $x^{(i)}=x^{(1)}$ for all $i \in [1, d],$ the  inequality is actually an equality.
Set $S_{2,n}= S_{n}-M_1$.  By (\ref{fsn02}), we have for any $p \in [1, \infty]$ and any $x>0$,
\begin{align}
  {\mathbb P}\Big(\| S_n \|_p \geq    d^{1/p} x \Big)
   &\leq {\mathbb P}\Big( \| M_1\|_p \geq    d^{1/p}x/2\Big)+ {\mathbb P}\Big( \| S_{2,n}\|_p \geq    d^{1/p}x/2\Big) \nonumber \\
   &\leq {\mathbb P}\Big( \| M_1\|_\infty \geq     x/2\Big)+ {\mathbb P}\Big( \| S_{2,n}\|_\infty \geq     x/2\Big) \nonumber
   \\
  &  \leq \sum_{i=1}^d \mathbb{P} \Big( | M_1^{(i)}  |  \geq    x/2 \Big ) +  \sum_{i=1}^d\mathbb{P} \Big( | S_{2,n}^{(i)}  |  \geq    x/2 \Big ) \nonumber \\
   & \leq d \max_{1\leq i \leq d} \mathbb{P} \Big( | M_1^{(i)}  |  \geq    x/2 \Big ) +   d \max_{1\leq i \leq d}  \mathbb{P} \Big( | S_{2,n}^{(i)}  |  \geq    x/2 \Big ). \label{ineq3.3}
\end{align}
Using the inequality
\begin{eqnarray}\label{ines4562}
|M_1^{(i)}| \leq \|M_1 \|_p \leq K_{ 1,n} G_{X_1}(X_1)
\end{eqnarray}
 (cf.\ the point  2 of Lemma \ref{McD}), we have for any $p \in [1, \infty]$ and any $x>0$,  
\begin{align}
  {\mathbb P}\Big(\| S_n \|_p \geq   d^{1/p} x \Big)
    &\leq d \, {\mathbb P}\bigg(  G_{X_1}(X_1) \geq  \frac{  x }{  2  K_{ 1,n}} \bigg) + d \max_{1\leq i \leq d}  \mathbb{P} \Big( | S_{2,n}^{(i)}  |  \geq  \frac  x 2 \, \Big ) \nonumber \\
    &  =: I_1(   x)+ I_2(   x) \, .\label{f29}
\end{align}
Hence,  to dominate ${\mathbb P}\big(\| S_n \|_p \geq   d^{1/p} x \big)$,
we only need to establish  deviation
inequalities for $I_1(  x)$ and $I_2(   x)$.
The term $I_1(  x)$  represents  the direct influence of the initial distribution of the chain, and it
  will be most of the time  negligible.  For instance, when  $X_1=x_1$ is a deterministic point, we have $G_{X_1}(X_1)=0$ and $I_1( x)=0$ for any $x>0$.
The main difficulty is to give an upper bound for $I_2( x)$, which is the purpose of  the remaining of this section.
By the point 3 of Lemma \ref{McD},  the martingale differences  $(M_k)_{k \in [2, n]}$ satisfy  for all $i \in [1, d]$,
\begin{eqnarray}\label{ines532}
  |M_k^{(i)}| \leq \|M_k \|_p\leq  K_{k, n} [ \tau_k  G_\varepsilon(\varepsilon_k) + \xi_k ] ,\, \ k \in [2,n].
\end{eqnarray}
Notice that $ S_{2,n}^{(i)} =\sum_{k=2}^n M_k^{(i)}.$
Since that   $G_\varepsilon(\varepsilon_k), k \in [2,n], $ are  i.i.d.\,random variables,
(\ref{ines532})  plays  an important role for estimating $\mathbb{P} \big( | S_{2,n}^{(i)}  |  \geq    x/2 \big )$
and  thus  $I_2( x)$. Assuming various moment conditions on  $G_\varepsilon(\varepsilon)$, we can obtain different upper bounds for $I_2(  x)$.

In the sequel, denote   by  $c_{p,d}$ and $c_{p,d}'$ positive constants,  which may depend  on the constants $ p, d,$   $\alpha, \rho$  and $\eta$ but do not depend on $x$ and $n$.

\subsection{Bernstein type bound}\label{Bersec}
In this subsection, we are interested in establishing a deviation inequality for $S_n$ under the Bernstein condition.
We refer to  de la Pe\~{n}a \cite{D99} for related inequalities:
in this paper  similar tight Bernstein type inequalities for  martingales are proved.
Using Lemma \ref{McD}, we get the following proposition.
\begin{prop}\label{prof3.1} Let $p \in [1, \infty].$
Assume that   there exist two positive constants   $H_1$ and $ A_1$  such that  for any integer $k\geq 2$,
\begin{equation}\label{Bernsteinmoment}
 {\mathbb E} \big[  \big(  G_\varepsilon(\varepsilon)\big)^k\big] \leq
 \frac {k!}{2}   H_1^{k-2} A_1 \, .
\end{equation}
Denote
$$V_n^2 = (1+ A_1)  \sum_{k=2}^n \big(2K_{ k,n}(\tau_{k}+\xi_k)\big)^2 \ \ \ \
 \textrm{and}  \ \ \ \ \delta_n=\max\{  2K_{ k,n}(\tau_{ k} H_1+\xi_k), k=2,\ldots,n\}.$$
Then  %for any $t \in [0, \delta_n^{-1})$ and $i \in [1, d]$,
%\begin{equation}\label{maindfs}
%  \mathbb{E}\,[e^ { \pm t  S_n^{(i)}   } ]\leq  \exp \left (\frac{t^2 V }{2 (1- t\,\delta_n )} \right )\, .
%\end{equation}
%Consequently,
for any $x> 0$,
\begin{eqnarray}
  {\mathbb P}\Big(\| S_n \|_p \geq  d^{1/p} x\Big)
  &\leq&I_1(x)+ 2d \exp \left\{- \frac{(x/2)^2}{V_n^2 (1+\sqrt{1+ x \delta_n/V_n^2 })+\delta_n x /2 }\right\}\, \label{Berie1} \\
 &\leq&I_1(x)+ 2d  \exp \left\{ -\frac{(x/2)^2}{2 \, V_n^2  + \delta_n x   }\right\}\, . \label{Berie2}
\end{eqnarray}
Assume moreover that there exits a positive constant $c$ such that for any $x>0,$
\begin{eqnarray}\label{ine35}
{\mathbb P}\big(  G_{X_1}(X_1) \geq x\big) \leq c^{-1} \, e^{-c\, x}.
\end{eqnarray}
Then inequality (\ref{Berie2}) implies that:
\begin{description}
  \item[\textbf{[i]}]  If (\ref{label1.5})   is satisfied,  then  for any $x>0,$
\begin{eqnarray}\label{th02s}
\limsup_{n\rightarrow \infty}\frac{1}{n} \ln \mathbb{P}\left(\| S_n \|_p \geq n x   \right)
    \ \leq \    - \, c_{p,d} \, \big( x \mathbf{1}_{\{x\geq 1\}} + x^2 \mathbf{1}_{\{0< x  <  1\}} \big) .
\end{eqnarray}
  \item[\textbf{[ii]}]  If  (\ref{label1.6}) is satisfied with $\alpha \in (0, 1/2)$, then  for any $x>0,$
\begin{eqnarray}\label{thd03s}
  \limsup_{n\rightarrow \infty}\frac{1}{n^{1-2\alpha} } \ln \mathbb{P}\left(\| S_n \|_p \geq n x   \right)
   \ \leq \  - c_{p,d} \, x^2 .
\end{eqnarray}
\item[\textbf{[iii]}]  If   (\ref{label1.7}) is satisfied,  then  for any $x>0,$
\begin{eqnarray}\label{th04s}
\limsup_{n\rightarrow \infty}\frac{1}{n} \ln \mathbb{P}\left(\| S_n \|_p \geq n x   \right)
    \ \leq \    - \, c_{p,d} \,  x   .
\end{eqnarray}
\end{description}
\end{prop}

If either (\ref{label1.5}) or (\ref{label1.7}) holds, then from   (\ref{th02s}) and (\ref{th04s}), it is easy to see that
$S_n$ admits the classical large deviation convergence rate $e^{-nc_x }$, where $c_x>0$.  Moreover, as $x > x^2$ for $0<x<1,$
the large deviation convergence rate in
(\ref{th04s}) is  better than that in  (\ref{th02s}).
 Under the condition (\ref{label1.6}),
the large deviation convergence rate for $S_n$   becomes  much worse, as shown by \eref{thd03s}.

As mentioned above, when $X_1 = x_1$ a.s. is deterministic,  it follows that $G_{X_1} (X_1) = 0$ a.s., and so
the condition (\ref{ine35}) holds for any constant $c \in (0, 1]$.

%From (\ref{Berie2}), we have the following result for confidence regions of $\frac{1}{n} S_n$.  For any given $\varepsilon>0$ and $\delta \in (0, 1)$,
%by (\ref{Berie2}),  we deduce that when $$n \geq  \frac{d^{1/p}\delta_n}{\varepsilon} \Big(\ln \frac{2d}{\varepsilon} \Big) + \frac{d^{1/p}}{\varepsilon}\sqrt{      \delta_n^2 \Big(\ln \frac{2d}{\delta} \Big)^2+    2V_n \Big(\ln\frac{2d}{\delta} \Big)\ },$$
%it holds  $\big\| \frac{1}{n} S_n \big \|_p< \varepsilon $
%with probability at least $1-\delta.$
%
 
Since
$G_\varepsilon(y) \leq  \delta(y, y_0) + {\mathbb E}[\delta(y_0, \varepsilon )]$, it follows that
$
{\mathbb E} \big[ \big( G_\varepsilon(\varepsilon)\big)^k\big] \leq 2^k {\mathbb E} \big[ \big(  \delta(\varepsilon, y_0)\big)^k\big]  .
$
Hence, the following condition
\begin{equation}\label{Berbis}
{\mathbb E} \big[ \big(  \delta(\varepsilon, y_0)\big)^k\big] \leq
\frac {k!}{2} A(y_0)^{k-2} B(y_0), \ \ \ k\geq 2,
\end{equation}
implies  the condition   \eref{Bernsteinmoment} with $H_1=2A(y_0)$ and $ A_1=4B(y_0)$.

 In  Examples 3, 4 and 5,   when $\mathcal{X}=\mathcal{Y} = \mathbb{R}^d$,  we can take $ f(x_1, x_2, \ldots, x_n)=   \sum_{i=1}^{n}x_i, \| \cdot \|=\| \cdot \|_p $
and   $d(x , x')=\delta(x , x' )=\| x -x' \|_p.$   Then the condition \eref{Bernsteinmoment} is satisfied, provided that
\begin{equation} \label{cond21}
{\mathbb E}|\varepsilon^{(i)}|^k \leq
 \frac {k!}{2}   H_1^{k-2} A_1, \ \ \ i \in [1, d]\ \textrm{ and }\, k\geq 2.
 \end{equation}
 To show this, by (\ref{Berbis}) with $y_0=0$, we only need to prove that
\begin{equation} \label{condsa}
{\mathbb E} \|\varepsilon \|_p ^k  \leq
\frac {k!}{2}   (  H_1d^{(p+1)/p}  )^{k-2}   A_1 d^{2(p+1)/p},\ \ \ k \geq2.
\end{equation}
Clearly, it holds
$$\| x \|_p^k \leq \bigg(d^{1/p}\sum_{i=1}^d |x^{(i)}| \bigg)^k   \leq d^{k/p} d^{k-1}  \sum_{i=1}^d |x^{(i)}|^k,\ \ \ \ k\geq 2.  $$
Hence, by (\ref{cond21}), we have
\begin{equation} \nonumber
{\mathbb E} \|\varepsilon \|_p ^k  \leq d^{k/p} d^{k-1}  \sum_{i=1}^d  {\mathbb E}| \varepsilon^{(i)}|^k \leq
\frac {k!}{2}   (  H_1d^{(p+1)/p}  )^{k-2}   A_1 d^{2(p+1)/p},\ \ \ k \geq2,
\end{equation}
which gives (\ref{condsa}).

\noindent\emph{Proof.}
  Notice that $S_{2,n}^{(i)} = \sum_{k=2}^nM_k^{(i)}$ is a sum of martingale differences.
By  (\ref{ines532}),  we have  for any  $k \in [2, n],$
\begin{eqnarray*}
|M_k^{(i)}|^{j} \leq   K_{k, n}^j [ \tau_k  G_\varepsilon(\varepsilon_k) + \xi_k ]^j \leq  2^{j-1} K_{k, n}^j [ \tau_k^j  (G_\varepsilon(\varepsilon_k))^j + \xi_k^j ]
\leq   2^{j} K_{k, n}^j [ \tau_k^j  (G_\varepsilon(\varepsilon_k))^j + \xi_k^j ]
\end{eqnarray*}
 and so for any $t >0$,
$$
\mathbb{E}\,[e^{t   M_k^{(i)} } ] \leq 1+ \sum_{j=2}^{\infty} \frac{t^j}{j!}\, \mathbb{E}\,[   | M_k^{(i)} | ^j ]
\leq 1+ \sum_{j=2}^{\infty} \frac{t^j}{j!}\, (2K_{k, n})^j \left[ \tau_k^j \mathbb{E}\,\big[ \big(  G_\varepsilon(\varepsilon_k)\big)^j \big] + \xi_k^j\right].
$$
Using the condition (\ref{Bernsteinmoment}), we deduce that  for any  $k \in [2, n]$ and  any $t \in [0, \delta_n^{-1})$,
\begin{eqnarray}
\mathbb{E}\,[e^{t   M_k^{(i)} } ]
&\leq& 1+ \sum_{j=2}^{\infty} \frac{t^j}{j!}\, (2K_{k, n})^j \left[ \tau_k^j \frac{j!}{2}H_{1}^{j-2}A_1 + \xi_k^j\right] \nonumber
 \\
&=&   1+ \frac{1}{2}\sum_{j=2}^{\infty} t^j    (2K_{k, n})^j \left[\tau_k^jH_{1}^{j-2}A_1 + \frac{2}{j !}  \xi_k^j\right]   \nonumber
 \\
&\leq& 1+ \frac{1}{2} t^2(2K_{k, n})^2 \sum_{j=2}^{\infty}  t^{j-2} \, (2K_{k, n})^{j-2}  (A_1+1)  [\tau_k H_{1}  + \xi_k ]^{j-2}(\tau_k    + \xi_k)^2 \nonumber
\\
&\leq&    1+ \frac{t^2 (A_1 +1) \big(   2 K_{k, n} (\tau_k +\xi_k)\big)^2  }{2 (1  -t\, \delta_n )}. \label{fin35}
\end{eqnarray}
Applying the inequality $1+u \leq e^u $  for $u\ge0$ to \eqref{fin35},  we deduce that  for any  $k \in [2, n]$ and  any $t \in [0, \delta_n^{-1})$,
\begin{eqnarray}
\mathbb{E}\,[e^{t M_k^{(i)} }|\mathcal{F}_{k-1}]  \leq  \exp \left\{  \frac{t^2  (1+A_1)   \big(   2 K_{k, n} (\tau_k +\xi_k)\big)^2  }{2 (1  -t\, \delta_n )}  \right\}.
\end{eqnarray}
By the tower property of conditional expectation and the last  inequality, it is easy to see that  for any $n\geq2$ and any $t \in [0, \delta_n^{-1})$,
\begin{eqnarray}
\mathbb{E}\,\big[e^{ tS_{2,n}^{(i)} } \big] &=&  \mathbb{E}\,\big[ \mathbb{E}\, [e^{ tS_{2,n}^{(i)} } |\mathcal{F}_{n-1}  ] \big]\\
\nonumber
& =&   \mathbb{E}\,\big[ e^ { tS_{2, n-1}^{(i)} } \mathbb{E}\,  [e^ { tM_{n}^{(i)} }  |\mathcal{F}_{n-1}  ] \big] \nonumber  \\
&\leq & \mathbb{E}\,\big[ e^ { tS_{2, n-1}^{(i)} } \big] \exp \left\{ \frac{t^2  (1+A_1)  \big(   2 K_{k, n} (\tau_k +\xi_k)\big)^2  }{2 (1  -t\, \delta_n )}  \right\} \nonumber \\
&\leq &  \exp \left\{\frac{t^2 V_n^2 }{2 (1- t\,\delta_n)} \right\}. \label{fsdfs}
\end{eqnarray}
Clearly, the same bound holds for $\mathbb{E}\,\big[e^{ -tS_{2,n}^{(i)} } \big]$. Returning to \eqref{f29},  by the exponential Markov  inequality,  we have
for any $x> 0$ and any $t \in [0, \delta_n^{-1})$,
\begin{eqnarray}
I_2(x)  &\leq&d \max_{1\leq i \leq d} \mathbb{E}\,  \Big[ \exp \big\{ tS_{2,n}^{(i)} -  \frac1 2 t x )  \big\} +  \exp \big\{ -tS_{2,n}^{(i)} - \frac1 2 tx    \big\}  \Big] \nonumber \\
 &  \leq&   2d   \exp \left\{-\frac1 2 t\,x  +  \frac{t^2 V_n^2 }{2 (1- t\,\delta_n)}    \right\}\, . \label{fines}
\end{eqnarray}
The last bound reaches its  minimum  at $$t=t(x):= \frac{ x/V_n^2 }{  x\delta_n/V_n^2  +1 + \sqrt{1+ x\delta_n/V_n^2 } }\, .$$ Substituting $t=t(x)$ in (\ref{fines}),  we obtain  for any $x> 0,$
\begin{eqnarray}
I_2(x)  &\leq&  2d   \exp \left\{- \frac{(x/2)^2}{V_n^2 (1+\sqrt{1+ x \delta_n/V_n^2 })+x  \delta_n/2 }\right\}\, \nonumber\\
 &\leq&  2d   \exp \left\{ -\frac{(x/2)^2}{2  V_n^2  +x  \delta_n   }\right\}\, ,\nonumber
\end{eqnarray}
where the last line follows by the inequality $\sqrt{1+ x\,\delta_n/V_n^2 }  \leq 1+x\, \delta_n/(2V_n^2) $.
Applying the upper bounds to (\ref{f29}), we obtain the inequalities (\ref{Berie1}) and (\ref{Berie2}).
\medskip

Condition (\ref{ine35}) implies that
\begin{eqnarray}
I_1 (x)  \leq  d c^{-1} \, \exp \bigg\{- \frac{c\,   x }{  2  K_{ 1,n} } \bigg\}.\label{fgjsfds1}
\end{eqnarray}
  If the  condition (\ref{label1.5}) is satisfied,  then we  have $K_{ 1,n}=O(1)$, $V_n^2=O(n)$ and $\delta_n=O(1)$ as $n\rightarrow \infty$,  by the point 4 of Lemma \ref{McD}.
 Applying (\ref{fgjsfds1}) to (\ref{Berie2}), we deduce that for any $x>0,$
\begin{eqnarray*}
  \mathbb{P}\left(\| S_n \|_p \geq n x   \right)
    \leq  d c^{-1} \, \exp \bigg\{-  c'_{p,d} nx \bigg\}+ 2d   \exp \left\{ - \,  \frac{ c'_{p,d} (n x  )^2}{  c_{p,d} ( n  +   n x )   }\right\}.
\end{eqnarray*}
 This last inequalitiy implies  (\ref{th02s}).
\medskip

  If the  condition (\ref{label1.6}) is satisfied, from the  point 5 of Lemma \ref{McD}, then we have $K_{ 1,n}=O(1)$,
   $$V_n^2=O(1)\sum_{k=2}^n  k^{2\alpha}=O(1)n^{1+2\alpha}\sum_{k=2}^{n} (\frac{k}{n})^{2\alpha}\frac{1}{n}=  O(n^{1+2\alpha})$$ and $\delta_n=O(n^{ \alpha})$ as $n\rightarrow \infty.$
 Applying (\ref{fgjsfds1}) to (\ref{Berie2}), we deduce that for any $\alpha \in (0, 1/2)$ and  any $x>0,$
\begin{eqnarray*}
  \mathbb{P}\left(\| S_n \|_p \geq n x   \right)
    \leq  d c^{-1} \, \exp \bigg\{-  c'_{p,d} nx \bigg\}+   2d  \exp \left\{ -\, \frac{c'_{p,d} (nx )^2}{  c_{p,d} ( n^{1+2\alpha}   + n x   n^{ \alpha} ) }\right\} .
\end{eqnarray*}
  The last inequality implies  (\ref{thd03s}).\medskip

 If the condition (\ref{label1.7}) is satisfied,  by the point 6 of Lemma \ref{McD}, then we have $K_{ 1,n}=O(1)$,
\begin{eqnarray}\label{fsdf36}
 V_n^2=O(1)\sum_{k=1}^n   \frac{1}{k^{2\alpha}}= \left\{ \begin{array}{ll}
  O(n^{1-2\alpha})  & \textrm{if $0\leq \alpha < 1/2$ }\\
  O(\ln n)  & \textrm{if $ \alpha = 1/2$ }\\
O(1)  & \textrm{if $ 1/2 < \alpha \leq1$ }
\end{array} \right.
\end{eqnarray}
 and  $\delta_n=O(1)$ as $n\rightarrow \infty.$
 Applying (\ref{fgjsfds1}) to (\ref{Berie2}), we deduce that for any $x>0,$
\begin{eqnarray*}
  \mathbb{P}\left(\| S_n \|_p \geq n x   \right)
    \leq  d c^{-1} \, \exp \bigg\{-  c'_{p,d} nx \bigg\}+  2d   \exp \left\{- \, \frac{ c'_{p,d} (n x )^2}{  c_{p,d} (  n^{1-2\alpha}\vee \ln n   +   n x )     }\right\}.
\end{eqnarray*}
 The last line implies  (\ref{th04s}). This completes the proof of Proposition \ref{prof3.1}.
  \hfill\qed

\subsection{Semi-exponential bound}
If both $G_{X_1}(X_1)$ and $G_{\varepsilon}(\varepsilon)$
have semi-exponential moments, the following proposition holds.  It  can be compared to the
corresponding results in Borovkov \cite{Bor} for partial sums of independent random variables,  Merlev\`{e}de {\it et al.}\ \cite{RPR10} for
partial sums of weakly
dependent sequences, Lesigne and  Voln\'{y} \cite{LV01} and Fan {\it et al.}\ \cite{FGL17} for martingales.
\begin{prop}\label{findsa}  Let $p \in [1, \infty] $ and $q\in (0, 1)$.
Assume that there exists  a  positive constant  $A_1$   such that
\begin{equation}\label{laplace2}
  {\mathbb E} \big[ \big( G_\varepsilon(\varepsilon)\big)^2\exp  \big\{  \big( G_\varepsilon(\varepsilon)\big)^q\big\} \big] \leq A_1
 \, .
\end{equation}
Denote
$$V_n^2 =  2e \sum_{k=2}^n  K_{k,n}  ^2 \Big( \tau_{k}^2  A_1   + \xi_k^2 \mathbb{E}[e^{| G_\varepsilon(\varepsilon)  |^q   }]   \Big)   \ \ \ \   \textrm{and}  \ \  \ \
 \delta_n=\max\Big\{    K_{ k,n}\tau_{ k},  K_{ k,n}  \xi_k, k=2,\ldots,n \Big\}.$$
If $V_n \geq 1 ,$ then  for any $x>0,$
\begin{eqnarray}\label{sens2}
  \mathbb{P}\left(\| S_n \|_p \geq  d^{1/p} x   \right)
 \!\! &\leq& \!\! I_1(x)+  4d  \exp\bigg\{- \, \frac{ (x/2)^2 }{2(V_n^2 + (x/2)^{2-q}\delta_n^{q}) }  \bigg\} .
\end{eqnarray}
Assume moreover that there exits a positive constant $c$ such that for any $x>0,$
\begin{eqnarray}\label{ifnelh5}
{\mathbb P}\big(  G_{X_1}(X_1) \geq x\big) \leq c^{-1} \, e^{-c\, x^q}.
\end{eqnarray}
Then inequality (\ref{sens2}) implies that:
\begin{description}
  \item[\textbf{[i]}]  If (\ref{label1.5}) or (\ref{label1.7}) is satisfied,  then for any $x>0,$
\begin{eqnarray}\label{th0ds3s}
  \limsup_{n\rightarrow \infty}\frac{1}{n^{q} } \ln \mathbb{P}\left(\| S_n \|_p \geq n x   \right)
    \ \leq \   -c_{p,d}  x^q   .
\end{eqnarray}
  \item[\textbf{[ii]}]  Assume that (\ref{label1.6}) is satisfied with $\alpha \in (0, 1/2)$. If $0< \alpha < \frac{1-q}{2-q}$,  then for any $x>0,$
  \begin{displaymath}
 \limsup_{n\rightarrow \infty}\frac{1}{n^{q(1-\alpha)}} \ln \mathbb{P}\left(\| S_n \|_p \geq n x   \right)
    \leq  -c_{p,d}  x^q.
\end{displaymath}
If $  \alpha= \frac{1-q}{2-q}$,  then for any $x>0,$
  \begin{displaymath}
 \limsup_{n\rightarrow \infty}\frac{1}{n^{q/(2-q)}} \ln \mathbb{P}\left(\| S_n \|_p \geq n x   \right)
    \leq  -c_{p,d} \big( x^q \mathbf{1}_{\{x\geq 1\}} + x^2 \mathbf{1}_{\{0< x  < 1\}} \big).
\end{displaymath}
If $ \frac{1-q}{2-q}< \alpha < \frac1 2$, then for any $x>0,$
\begin{displaymath}
 \limsup_{n\rightarrow \infty}\frac{1}{n^{1-2\alpha} } \ln \mathbb{P}\left(\| S_n \|_p \geq n x   \right)
    \leq    -c_{p,d}  x^2.
\end{displaymath}
\end{description}
\end{prop}

If either (\ref{label1.5}) or (\ref{label1.7}) is satisfied,  from (\ref{th0ds3s}), it is easy to see that
the large deviation convergence rate is  the same as  the classical one. On the contrary under the condition (\ref{label1.6}),
this convergence rate becomes worsen as $\alpha$ increases.
\medskip

\noindent\emph{Proof.}
Notice that $S_{2,n}^{(i)}/\delta_n = \sum_{k=2}^nM_k^{(i)}/\delta_n $ is a sum  of martingale differences. By  (\ref{ines532}) and the condition (\ref{laplace2}), it is  easy to see that  for any $k \in [2, n],$
\begin{eqnarray*}
\mathbb{E}\,[ (M_k^{(i)}/\delta_n)^2e^{|M_k^{(i)}/\delta_n|^q } ] &\leq&  \delta_n^{-2}  \mathbb{E}\,[ (K_{k,n}[\tau_{k}G_\varepsilon(\varepsilon_k) +\xi_k])^2e^{|K_{k,n}[\tau_{k} G_\varepsilon(\varepsilon_k) + \xi_k ] /\delta_n |^q } ]  \nonumber\\
&\leq& 2(K_{k,n}\delta_n^{-1})^2 \mathbb{E}\,[ ((\tau_{k}G_\varepsilon(\varepsilon))^2 + \xi_k^2 ) e^{| G_\varepsilon(\varepsilon)  |^q +1 } ]  \nonumber\\
&\leq& 2e (K_{k,n}\delta_n^{-1})^2 \Big( \tau_{k}^2 \mathbb{E}\,[  (G_\varepsilon(\varepsilon))^2 e^{| G_\varepsilon(\varepsilon)  |^q   }  ]   + \xi_k^2 \mathbb{E}[e^{| G_\varepsilon(\varepsilon)  |^q   }]   \Big)  \nonumber\\
&\leq&  2e (K_{k,n}\delta_n^{-1})^2 \Big( \tau_{k}^2  A_1   + \xi_k^2 \mathbb{E}[e^{| G_\varepsilon(\varepsilon)  |^q   }]   \Big)  .
\end{eqnarray*}
If $V_n \delta_n^{-1} \geq 1,$ using inequality (2.7) of Fan \emph{et al.}\ \cite{FGL17}, then  we have for any $t>0,$
\begin{eqnarray}  \label{fdisaxs}
  \mathbb{P}\left( | S_{2,n}^{(i)}/\delta_n  |  \geq  t \, \right)
 \!\! &\leq& \!\!   4   \exp\bigg\{- \, \frac{ t^2 }{2(V_n^2\delta_n^{-2} + t^{2-q}) }  \bigg\}.
\end{eqnarray}
Substituting $t= x/(2\delta_n) $ in (\ref{fdisaxs}),  we get for any $x>0,$
\begin{eqnarray*}
  \mathbb{P}\left( | S_{2,n}^{(i)}  |  \geq  x/2 \, \right)
 \!\! &\leq& \!\!   4  \exp\bigg\{- \, \frac{ (x/2)^2 }{2(V_n^2 + (x/2)^{2-q}\delta_n^{q}) }  \bigg\} .
\end{eqnarray*}
From (\ref{f29}) and the last inequality, we obtain  the desired inequality (\ref{sens2}).\medskip

Condition (\ref{ifnelh5}) implies that
\begin{eqnarray}
I_1 (x)  \leq  d c^{-1} \, \exp \bigg\{- c\, \Big(\frac{   x }{  2  K_{ 1,n} } \Big)^q \bigg\}.\label{fgjslds2}
\end{eqnarray}
  If  the condition (\ref{label1.5}) is satisfied, then by the point 4 of Lemma \ref{McD},  we obtain $K_{ 1,n}=O(1)$,  $V_n^2=O(n)$ and $\delta_n=O(1)$ as $n\rightarrow \infty.$
  Applying   inequality  (\ref{fgjslds2}) together with (\ref{sens2}), we deduce that for any $x>0,$
\begin{eqnarray*}
  \mathbb{P}\left(\| S_n \|_p \geq n x   \right)
    \leq  d c^{-1} \,  \exp \bigg\{- c'_{p,d} (nx)^q   \bigg\}+ 4d \exp\bigg\{- \, \frac{  c'_{p,d} (n x )^2 }{   c_{p,d} (  n + (n x )^{2-q} )  }  \bigg\}.
\end{eqnarray*}
From the last inequality, we get  (\ref{th0ds3s}).\medskip

  If  the condition (\ref{label1.6}) is satisfied,  then by the point 5 of Lemma \ref{McD},  we have $K_{ 1,n}=O(1)$,  $V_n^2=O(n^{1+2\alpha})$ and $\delta_n=O(n^{ \alpha})$ as $n\rightarrow \infty.$
  Applying (\ref{fgjslds2}) to inequality (\ref{sens2}), we deduce that for  any $\alpha \in (0, 1/2)$ and any $x>0,$
\begin{displaymath}
  \mathbb{P}\left(\| S_n \|_p \geq n x   \right)
    \leq d c^{-1} \,  \exp \bigg\{- c'_{p,d} (nx)^q   \bigg\}+ 4 d   \exp\bigg\{- \, \frac{ c'_{p,d} (n x )^2 }{    c_{p,d} (  n^{1+2\alpha} + (n x )^{2-q}  n^{ q\alpha} ) }  \bigg\}.
\end{displaymath}
  The last inequality implies the point [ii] of Proposition \ref{findsa}. Note  that when $0< \alpha <  \frac{1-q}{2-q}$, we have $n^{2-q+ q\alpha} > n^{1+2\alpha}$,
  while when $\frac{1-q}{2-q}< \alpha < \frac1 2$, then $n^{2-q+ q\alpha} < n^{1+2\alpha}$.

 If   the condition (\ref{label1.7}) is satisfied,  then by the point 6 of Lemma \ref{McD},  we have $K_{ 1,n}=O(1)$, \eref{fsdf36}
 and  $\delta_n=O(1)$ as $n\rightarrow \infty.$
Applying again inequalities (\ref{fgjslds2}) and (\ref{sens2}), we deduce that for any $x>0,$
\begin{eqnarray*}
  \mathbb{P}\left(\| S_n \|_p \geq n x   \right)
    \leq d c^{-1} \, \exp \bigg\{- c'_{p,d} (nx)^q   \bigg\}+ 4 d   \exp\bigg\{- \, \frac{ c'_{p,d} (n x )^2 }{  c_{p,d} (    n^{1-2\alpha}\vee \ln n  + (n x )^{2-q} ) }  \bigg\}.
\end{eqnarray*}
Inequality  \eqref{th0ds3s}  is an easy consequence of this last inequality.
  \hfill\qed

\subsection{Fuk-Nagaev type bound}\label{Fuksec}
If the martingale differences $(M_i)_{i\geq2}$ admit finite $q$th order moments ($q\geq 2$), then we have the following Fuk-Nagaev type inequality (cf.\ Corollary $3'$ of Fuk \cite{F73} and Nagaev \cite{N79}).
\begin{prop}\label{mineq}
 Let  $p\in[1,\infty]$ and   $q\in [2, \infty)$.  Assume that there exist two positive constants   $A_1 $
and $B_1(q)$   such that
\begin{eqnarray}
{\mathbb E}\big[\big(G_\varepsilon(\varepsilon)\big)^2 \big] \leq A_1\,
 \,
\quad   \text{and}  \quad \ \, \ \,
 {\mathbb E}\big[\big(G_\varepsilon(\varepsilon)\big)^q \big] \leq B_1(q)
  \,  . \label{sdvcx}
\end{eqnarray}
Denote
$$
V_n^2 =  2  \sum_{k=2}^n K_{k, n}^2 \Big( \tau_k ^2 A_1  +  \xi_k^2\Big) \ \ \  \ \ \ \textrm{and}  \ \ \ \ \
   H_n(q)=  2^{q-1} \sum_{k=2}^n K_{k, n}^q \Big( \tau_k ^q B_1(q) +  \xi_k^q\Big) .
$$
Then for any $x  > 0$,
\begin{eqnarray}\label{fuki}
{\mathbb P}\Big(\| S_n \|_p \geq  d^{1/p} x\Big)   \leq I_1(x) + 2^{q+1}d\Big(1+ \frac2 q \Big)^q \frac{H_n(q)}{x^q}   + 2d \exp\left\{ - \frac{ x^2}{2(q+2)^2 e^q V_n^2}   \right\}.
\end{eqnarray}
Assume moreover that there exists a positive constant $c$ such that for any $x>0,$
\begin{eqnarray}\label{ifgnelh5}
{\mathbb P}\big(  G_{X_1}(X_1) \geq x\big) \leq c  \,   x^{ - q} ,
\end{eqnarray}
 then  inequality (\ref{fuki}) implies that:
\begin{description}
  \item[\textbf{[i]}]  If (\ref{label1.5})  is satisfied,  then for any $x>0,$
\begin{eqnarray}\label{ineq01}
  \mathbb{P}\left(\| S_n \|_p \geq n x   \right)
    \leq   \frac{c_{p, d}}{  x^q  }\cdot \frac{1}{  n^{q-1}}.
\end{eqnarray}
  \item[\textbf{[ii]}]  If  (\ref{label1.6}) is satisfied with $0< \alpha <  \frac{1}{2} $, then for any $x>0,$
\begin{eqnarray}\label{ineq02}
  \mathbb{P}\left(\| S_n \|_p \geq n x   \right)
    \leq    \frac{c_{p, d}}{  x^q  }\cdot  \frac{1}{\   n^{q-1-\alpha q}} .
\end{eqnarray}
\item[\textbf{[iii]}]  If  (\ref{label1.7}) is satisfied, then for any $x>0,$
\begin{eqnarray}\label{ineq03}
 \mathbb{P}\left(\| S_n \|_p \geq n x   \right)
    \leq   \left\{ \begin{array}{ll}
\displaystyle   \frac{c_{p, d}}{  x^q  }\cdot  \frac{1}{\   n^{q-1+\alpha q}} \,,  & \textrm{if $ 0<\alpha < \frac1 q,$ }\\
\\
\displaystyle   \frac{c_{p, d}}{  x^q  }\cdot  \frac{\ln n}{\   n^{q}}\,,  & \textrm{if $ \alpha = \frac1 q,$ }\\
\\
\displaystyle   \frac{c_{p, d}}{  x^q  }\cdot  \frac{1}{\   n^{q}}\,,  & \textrm{if $\frac1 q<  \alpha <1$.}
\end{array} \right.
\end{eqnarray}
\end{description}
\end{prop}

Under the condition (\ref{label1.5}), then from (\ref{ineq01}), it is easy to see that
the large deviation convergence rate is the same as  the classical one, which  is of order $n^{1-q}$ as $n \rightarrow \infty$.
When the condition (\ref{label1.6}) is satisfied, from (\ref{ineq02}), we find that the large deviation convergence rate becomes worse when $\alpha$ tends to $1/2$.
Moreover, the large deviation convergence rate is slower than the classical one.
Now, if the  condition (\ref{label1.7}) is satisfied, then the inequalities  (\ref{ineq03})  imply  that this convergence rate
is much better than the classical one.

\noindent\emph{Proof.}
By   \eref{ines532} and the  condition (\ref{sdvcx}), it follows that
\begin{eqnarray*}
\sum_{k=2}^n\mathbb{E}[|M_k^{(i)}|^q | \mathcal{F}_{k-1} ] &\leq&   \sum_{k=2}^n\mathbb{E}[ \big(K_{k, n}( \tau_k  G_\varepsilon(\varepsilon_k)  + \xi_k)\big)^q ]= 2^{q-1}\sum_{k=2}^n \Big(\big( K_{k, n} \tau_k \big)^q\mathbb{E}[( G_{\varepsilon}(\varepsilon ))^q ] + K_{k, n}^q \xi_k^q\Big)\\
&\leq& 2^{q-1} \sum_{k=2}^n K_{k, n}^q \Big( \tau_k ^q B_1(q) +  \xi_k^q\Big)
=   H_n(q).
\end{eqnarray*}
Notice that $H_n(2)=V_n^2$. Using the Corollary $3'$ in Fuk \cite{F73}, we have for any $x>0,$
\begin{eqnarray}
{\mathbb P}(  |S_{2,n}^{(i)}  |  \geq x/2 )   \leq  2^{q+ 1}\Big(1+ \frac2 q \Big)^q \cdot \frac{H_n(q)}{x^q}   + 2 \exp\left\{ - \frac{x^2}{2(q+2)^2 e^q V_n^2 }  \right\}.
\end{eqnarray}
Applying the last inequality to (\ref{f29}), we get the first desired inequality.

The condition (\ref{ifgnelh5}) implies that for any $x>0,$
\begin{eqnarray}\label{fsdhy3}
I_1 (x)  \leq c \, d     \Big(\frac{   x }{  2  K_{ 1,n} } \Big)^{-q}   .
\end{eqnarray}
  If the condition (\ref{label1.5}) is satisfied, then from the point 4 of Lemma \ref{McD},  we have $K_{1,n}=O(1)$, $V_n^2=O(n)$ and $H_n(q)=O(n)$ as $n\rightarrow \infty.$
Applying (\ref{fsdhy3}) to \eref{fuki}, we get for any $x  > 0$,
\begin{eqnarray*}
{\mathbb P}( \| S_n \|_p \geq  n x )   \leq  \frac{ c_{p,d}}{ \big(n x  \big)^{ q}}   +  \frac{c_{p,d}\,n}{(n x )^q}   + 2d \exp\left\{ -  c'_{p,d} \frac{(n x )^2}{ \, n}  \right\}.
\end{eqnarray*}
The inequality \eref{ineq01} follows from the last inequality.

  If the condition (\ref{label1.6}) is satisfied,  by the  point 5 of Lemma \ref{McD}, then we have $K_{1,n}=O(1)$, $V_n^2=O(n^{1+2\alpha})$ and $H_n(q)= O(1) \sum_{k=2}^n k^{\alpha q}= O(n^{1+\alpha q})$ as $n\rightarrow \infty.$
Applying (\ref{fsdhy3}) to \eref{fuki}, we get  for any $x  > 0$,
\begin{eqnarray*}
{\mathbb P}( \| S_n \|_p \geq  n x )   \leq \frac{ c_{p,d}}{ \big(n x  \big)^{ q}}  +   \frac{c_{p,d} \,n^{1+ \alpha q}}{(n x )^q}   + 2d \exp\left\{ - c'_{p,d} \frac{(n x )^2}{ \, n^{1+2\alpha  }}  \right\}.
\end{eqnarray*}
The inequality \eref{ineq02} follows from this last inequality. Note that if $\alpha < 1/2,$ then  the third  term in the right hand side of the last inequality tends to $0$ as $n\rightarrow \infty$.

 If the condition (\ref{label1.7}) is satisfied, then  by the  point 6 of Lemma \ref{McD},  we have $K_{ 1,n}=O(1)$, \eref{fsdf36}
and
\begin{displaymath}
 H_n(q)=O(1) \sum_{k=1}^n k^{-\alpha q}=  \left\{ \begin{array}{ll}
  O(n^{1- \alpha q})\,,   & \textrm{if $ 0 \leq \alpha < 1/q,$ }\\
  O(\ln n)\,,  & \textrm{if $ \alpha = 1/q,$ }\\
O(1)\,,   & \textrm{if $ 1/ q < \alpha   <1$ }
\end{array} \right.
\end{displaymath}
 as $n\rightarrow \infty.$
Similarly, we  prove that the inequality (\ref{fuki}) implies \eref{ineq03}.
  \hfill\qed

\subsection{von Bahr-Esseen type   bound}\label{VBdEB}
If the dominating random variables $G_{X_1}(X_1)$ and $G_\varepsilon(\varepsilon )$ admit only a finite moment with order $q \in [1,2]$,
we have the following von Bahr-Esseen type deviation bound.
\begin{prop}\label{VBEI}
Let  $p\in[1,\infty]$ and  $q \in [1, 2]$. Assume that  there exists a positive constant   $A_1(q) $  such that
\begin{equation}\label{fnafdsa58}
{\mathbb E} \big[  \big(   G_\varepsilon(\varepsilon)\big)^q\big] \leq A_1(q) \, .
\end{equation}
Denote
$$V_n(q)= 2^{q-1} \bigg[  K_{2, n} ^q (\tau_2^q A_1(q)  + \xi_2^q)
         + 2^{2-q}
    \sum_{k=3}^n  K_{k, n} ^q    \big( \tau_k^q A_1(q)  + \xi_k^q \big)  \bigg] .$$
Then for any $x>0,$
\begin{eqnarray} \label{vb02}
  {\mathbb P}\Big(\| S_n \|_p \geq  d^{1/p} x\Big)  \leq I_{1}(x) +  2^q d  \frac{V_n(q)}{x^q}  .
\end{eqnarray}
Assume moreover (\ref{ifgnelh5}).
 Then  inequality (\ref{vb02}) implies that:
\begin{description}
  \item[\textbf{[i]}]  If (\ref{label1.5})  is satisfied,  then \eref{ineq01} holds.
  \item[\textbf{[ii]}]  If  (\ref{label1.6}) is satisfied with $0< \alpha < 1- \frac{1}{q}$, then  (\ref{ineq02}) holds.
%\begin{eqnarray}\label{f3.35}
%  \mathbb{P}\left(\| S_n \|_p \geq n x   \right)
%    \leq    \frac{c_{p,d}}{  x^q  }  \frac{1}{\   n^{q-1-\alpha q}}.
%\end{eqnarray}
\item[\textbf{[iii]}]  If  (\ref{label1.7}) is satisfied, then \eref{ineq03} holds.
\end{description}
\end{prop}
\begin{rem}
The constant  $2^{2-q}$  in $V_n(q)$ can be replaced by the more precise constant  $\tilde C_q$  described in Proposition 1.8 of Pinelis \cite{P10}.
\end{rem}
\noindent\emph{Proof.} Notice that $S_{2,n}^{(i)}  = \sum_{k=2}^nM_k^{(i)} $ is a sum  of martingale differences. Using a refinement of the von Bahr-Esseen inequality (cf.\ Proposition 1.8 of Pinelis \cite{P10}), we get for any $q \in [1, 2],$
$$
{\mathbb E}  |S_{2,n}^{(i)} |^q  \leq {\mathbb E}  |M_2^{(i)} |^q +    2^{2-q} \sum_{k=3}^n {\mathbb E}  |M_k^{(i)} |^q   \,.
$$
By \eref{ines532} and \eref{fnafdsa58}, we deduce that for any $q \in [1, 2],$
\begin{eqnarray}
{\mathbb E}  |S_{2,n}^{(i)} |^q   &\leq&   K_{2, n}  ^q   {\mathbb E} \big[  \big(   \tau_2  G_\varepsilon(\varepsilon) +\xi_2 \big)^q\big]    + 2^{2-q}
    \sum_{k=3}^n  K_{k, n}^q {\mathbb E} \big[  \big( \tau_k  G_\varepsilon(\varepsilon) + \xi_k \big)^q\big]\,  \nonumber  \\
    &\leq&    2^{q-1} \bigg[  K_{2, n}^q   {\mathbb E} \big[  \big(   \tau_2  G_\varepsilon(\varepsilon))^q + \xi_2 ^q\big]
         + 2^{2-q}
    \sum_{k=3}^n  K_{k, n} ^q {\mathbb E} \big[  \big( \tau_k  G_\varepsilon(\varepsilon))^q + \xi_k  ^q\big]\,  \bigg] \nonumber  \\
    &\leq&  \,  \, 2^{q-1} \bigg[  K_{2, n} ^q (\tau_2^q A_1(q)  + \xi_2^q)
         + 2^{2-q}
    \sum_{k=3}^n  K_{k, n} ^q    \big( \tau_k^q A_1(q)  + \xi_k^q \big)  \bigg] \nonumber  \\
    &  =&  V_n(q) .  \label{fsdsg45}
\end{eqnarray}
 By Markov's inequality, we get   for any $x  > 0$,
\begin{eqnarray}
{\mathbb P}(  |S_{2,n}^{(i)}  |  \geq x/2 )   &\leq&  2^q \frac{ {\mathbb E}|   S_{2,n}^{(i)}  |^q}{x^q}  \nonumber \\
&\leq& 2^q \frac{V_n(q)}{x^q}  .
\end{eqnarray}
 Applying the last inequality to (\ref{f29}), we get the desired inequality (\ref{vb02}).
The remaining of the proof is similar to the proof of Proposition \ref{mineq}.  \hfill\qed

Next, we consider the case where the random variables
$G_{X_1}(X_1)$ and $G_\varepsilon(\varepsilon)$ have only a weak moment.
Recall that for any real-valued random variable $Z$ and any $q\geq 1$,  the weak
moment of order $q$ is defined  by
\begin{equation}\label{weakp}
\|Z\|_{w,q}^q=\sup_{x>0} x^q{\mathbb P}( |Z | >x)\, .
\end{equation}
When  the variables $G_{X_1}(X_1)$ and
$G_\varepsilon(\varepsilon)$ have only a weak moment of order $q \in (1,2)$, we have the following deviation inequality.

\begin{prop}\label{weakmom}
Let  $p\in[1,\infty]$ and  $q \in (1, 2)$. Assume that  there exists a positive constant   $A_1(q) $  such that
\begin{equation}\label{weakVB}
  \big \|G_\varepsilon(\varepsilon)\big\|_{w,q}^q\leq A_1(q)
  \, .
\end{equation}
Then for any $x>0,$
\begin{equation}\label{weakVBE}
{\mathbb P}\Big(\| S_n \|_p \geq  d^{1/p} x\Big)  \leq   I_1(x) +  C_{d,q}\frac{ B(n,q)}{x^q}\, ,
\end{equation}
where
 $$
 C_{d,q}=2^{2+q}d\Big( \frac{ q  }{q-1} +\frac{2    }{2-q}\, \Big)
 \ \ \  \ \ \
\textrm{and}
   \ \ \  \ \ \
B(n,q)=  \sum_{k=2}^n (2K_{k, n}  )^q (   \tau_k ^q A_1(q)+ \xi_k ^q )   .
$$
Assume moreover $ \big \|G_{X_1}(X_1)\big\|_{w,q}^q  < \infty .$  Then inequality (\ref{weakVBE}) implies that:
\begin{description}
  \item[\textbf{[i]}]  If (\ref{label1.5})  is satisfied,  then \eref{ineq01} holds.
  \item[\textbf{[ii]}]  If  (\ref{label1.6}) is satisfied with $0< \alpha < 1- \frac{1}{q}$, then    (\ref{ineq02})  holds.
\item[\textbf{[iii]}]  If  (\ref{label1.7}) is satisfied, then \eref{ineq03} holds.
\end{description}
\end{prop}

\noindent\emph{Proof.} By Proposition 3.3 of Cuny,   Dedecker and  Merlev\`{e}de \cite{CDM17}, we have for any $x>0,$
$$
 {\mathbb P}\big( | S_{2,n}^{(i)}  | \geq   x/2 \big)  \leq    \frac{C_q}{x^q} \sum_{k=2}^{n} \| M_k^{(i)}\|_{w,q}^q  ,
$$
where $C_q=2^{2+q} ( \frac{ q}{q-1}  + \frac{2}{2-q}  ).$   By (\ref{f29}), we have for any $x>0,$
\begin{equation}\label{wedsBf}
{\mathbb P}\big(\| S_n \|_p \geq d^{1/p} x \big)  \leq   I_1(x) +    \frac{ C_{d,q}}{x^q} \max_{1\leq i \leq d}\sum_{k=2}^{n} \| M_k^{(i)}\|_{w,q}^q .
\end{equation}
Using (\ref{ines532}),  we have for any $k \in [2, n],$
\begin{eqnarray*}
\| M_k^{(i)}\|_{w,q}^q &\leq& \sup_{x>0} x^q{\mathbb P}(   K_{k, n} \tau_k  G_\varepsilon(\varepsilon) + K_{k, n} \xi_k  >x)   \\
&\leq&  (2K_{k, n} \xi_k)^q+ \sup_{x>  2K_{k, n}} x^q{\mathbb P}(   K_{k, n} \tau_k  G_\varepsilon(\varepsilon)  > x/2 ) \\
&\leq&  (2K_{k, n} \xi_k)^q + \sup_{x> 0} (2K_{k, n} \tau_kx)^q{\mathbb P}(      G_\varepsilon(\varepsilon)  > x  ) \\
&\leq&(2K_{k, n} \xi_k)^q +    (2K_{k, n} \tau_k   )^q  \big \|G_\varepsilon(\varepsilon)\big\|_{w,q}^q  \\
&\leq&  (2K_{k, n} \xi_k)^q +    (2K_{k, n} \tau_k   )^q A_1(q).
\end{eqnarray*}
Returning to (\ref{wedsBf}), we get for any $x>0,$
\[
{\mathbb P}\big(\| S_n \|_p \geq d^{1/p} x \big)  \leq   I_1(x) +     \frac{ C_{d,q} }{x^q}   \sum_{k=2}^n (2K_{k, n}  )^q (   \xi_k ^q +  \tau_k ^q A_1(q) )  ,
\]
which is exactly the first desired inequality.
 The remaining of the proof is similar to the proof of Proposition \ref{mineq}.  \hfill\qed

\subsection{McDiarmid type bound}
In this section, we consider the case where the increments $M_k$ are bounded in $L_{\infty}$-norm.
We shall make use of a refinement of the well-known McDiarmid inequality, which has been recently established by Rio \cite{R13}.
Following the notations in Rio \cite{R13},  denote
$$\ell(t)= (t - \ln t -1) +t (e^t-1)^{-1} + \ln(1-e^{-t}) \ \ \
\textrm{for all}\ \ t > 0, $$ and let
$$\ell^*(x)=\sup_{t>0}\big(xt- \ell (t) \big), \ \ \    x > 0, $$ be the Young transform of $\ell(t)$.
As quoted by Rio \cite{R13}, for any $x \in [0,1)$, it holds
\begin{equation}\label{is50}
\ell^*(x) \geq (x^2-2x) \ln(1-x) \geq 2x^2   .
\end{equation}
Denote by $ \varepsilon' $  an independent copy of
$ \varepsilon$, and $ X_1' $  an independent copy of
$ X_1.$

\begin{prop} \label{classic}
Let $p\in[1,\infty]$.
Assume that    there exists a positive constant  $T_1$ such that
\begin{equation}\label{in850}
  \|\delta(\varepsilon, \varepsilon' )\|_\infty \leq  T_1 .
\end{equation}
Let
$$
V_n^2=      \sum_{k=2}^n  K_{ k,n}^2(\tau_{ k} T_1 + \xi_k)^2
 \ \ \ \ \ \
\textrm{and}
  \ \ \ \ \ \
 D_n=    \sum_{k=2}^n K_{ k,n}(\tau_{ k} T_1 + \xi_k) .
$$
Then,  for any   $x \in [0, 2D_n]$,
\begin{eqnarray}\label{rioineq}
{\mathbb P}\Big(\| S_n \|_p \geq  d^{1/p} x\Big) \ \leq \ I_1(x)+ 2d \exp \left \{ -\frac{D^2_n}{V_n^2 } \ \ell^*\Big( \frac {x} {2D_n} \Big)\right\}  .
\end{eqnarray}
Consequently, for any $x \in [0, 2D_n]$,
\begin{eqnarray}
{\mathbb P}\Big(\| S_n \|_p \geq  d^{1/p} x\Big) & \leq  & I_1(x)+  2d \left ( \frac{D_n-x/2}{D_n}\right)^{ \frac{ D_n x -(x/2)^2}{V_n^2  }} \label{riosbou01} \\
&\leq &I_1(x)+ 2d\exp \left \{ -\frac{  x^2}{2V_n^2  }\label{riosbound}
 \right\}.
\end{eqnarray}
Assume moreover  $\|d(X_1, X'_1)\|_\infty< \infty.$ Then  we have:
\begin{description}
  \item[\textbf{[i]}]  If (\ref{label1.5})  is satisfied,  then for any $x>0,$
\begin{eqnarray}\label{ine010s}
 \limsup_{n\rightarrow \infty}\frac{1}{n} \ln \mathbb{P}\left(\| S_n \|_p \geq n x   \right)
    \leq    -c_{p,d} \, x^2 .
\end{eqnarray}
  \item[\textbf{[ii]}]  If  (\ref{label1.6}) is satisfied  with $\alpha \in (0, 1/2)$, then for any $x>0,$
\begin{eqnarray}\label{ine011s}
 \limsup_{n\rightarrow \infty}\frac{1}{n^{1-2\alpha}} \ln  \mathbb{P}\left(\| S_n \|_p \geq n x   \right)
    \leq    -c_{p,d} \, x^2 .
\end{eqnarray}
\item[\textbf{[iii]}]  Assume the condition (\ref{label1.7}).
If $0< \alpha < 1/2$, then for any $x>0,$
$$ \limsup_{n\rightarrow \infty}\frac{1}{n^{1+2\alpha}} \mathbb{P}\left(\| S_n \|_p \geq n x   \right)
    \leq  -c_{p,d} \, x^2\;.$$
If $\alpha = 1/2$, then for any $x>0,$
$$ \limsup_{n\rightarrow \infty}\frac{ \ln n}{n^2 } \mathbb{P}\left(\| S_n \|_p \geq n x   \right)
    \leq  -c_{p,d} \, x^2\;.$$
If $ 1/2 < \alpha  <1$, then for any $x>0,$
$$ \limsup_{n\rightarrow \infty}\frac{1}{n^2} \mathbb{P}\left(\| S_n \|_p \geq n x   \right)
    \leq  -c_{p,d} \, x^2\;.$$
\end{description}
\end{prop}

\noindent\emph{Proof.}
For any $k \in [2, n],$ let
$$
u_{k-1}^{(i)}(x_1, \ldots, x_{k-1})= \text{ess}\inf_{\varepsilon_k} g_k^{(i)}\big(x_1, \ldots, F_k(x_{k-1}, \varepsilon_k) \big),
$$
and
$$
v_{k-1}^{(i)}(x_1, \ldots, x_{k-1})= \text{ess}\sup_{\varepsilon_k} g_k^{(i)}\big(x_1, \ldots, F_k(x_{k-1}, \varepsilon_k) \big).
$$
From the proof of Proposition \ref{McD},  it follows that for any $k \in [2, n],$
$$
u_{k-1}^{(i)}(X_1, \ldots, X_{k-1}) \leq M^{(i)}_k \leq v_{k-1}^{(i)}(X_1, \ldots, X_{k-1})\, .
$$
By the proof of Proposition \ref{McD} and the  condition (\ref{in850}), we have
\[
v_{k-1}^{(i)}(X_1, \ldots, X_{k-1})-u_{k-1}^{(i)}(X_1, \ldots, X_{k-1}) \leq K_{k,n}(\tau_{k} T_1 + \xi_k) \, , \ \  k \in [2, n] .
\]
Now, with an argument similar to the proof of Theorem 3.1 of Rio \cite{R13} with    $\Delta_k = K_{ k,n}(\tau_{ k} T_1 + \xi_k), k \in [2, n]$,
we get
for any   $x \in [0,  2D_n]$,
\begin{eqnarray}
{\mathbb P}\big( | S_{2,n}^{(i)}|  \geq    x/2 \big) \ \leq \ 2  \exp \left \{ -\frac{D^2_n}{V_n^2  } \ \ell^*\Big( \frac {x} {2 D_n} \Big)\right\} \, .
\end{eqnarray}
Applying this last inequality to \eqref{f29}, we obtain (\ref{rioineq}).
 By the inequality
$\ell^*(x) \geq (x^2-2x) \ln(1-x), x \in [0,1)$, inequality (\ref{riosbou01})
follows from (\ref{rioineq}). Since for any $x \in [0,1)$, $(x^2-2x) \ln(1-x) \geq 2\,x^2$, it follows that for any $x \in [0,  2D_n]$,
\begin{eqnarray*}
\left ( \frac{D_n-x/2}{D_n}\right)^{ \frac{ D_n x -(x/2)^2}{V_n^2  }} \ \leq \   \exp \left \{ -\frac{x^2}{2 V_n^2  }
 \right\} ,
\end{eqnarray*}
which gives \eref{riosbound}.

 If $\|d(X_1, X'_1)\|_\infty< \infty,$ then we have for any $x>0,$
$$  I_1(x) \leq d \, {\mathbb P}\bigg( \|d(X_1, X'_1)\|_\infty \geq  \frac{  x }{  2  K_{ 1,n}} \bigg).$$
The last inequality implies that $I_1(x)=0$ for $x >2  K_{ 1,n}  \|d(X_1, X'_1)\|_\infty  .$

  If the condition (\ref{label1.5}) is satisfied, by the point 4 of Lemma \ref{McD},  then it holds $V_n^2=O(n)$.
Thus, inequality (\ref{riosbound}) implies that for any $x>0,$
\begin{eqnarray}
  \mathbb{P}\left(\| S_n \|_p \geq n x   \right)
    \leq I_1(nx d^{-1/p} )  +   2d\exp \left \{ - c_{p,d} \frac{ (x n )^2}{ \,  n }
 \right\} .
\end{eqnarray}
  The last inequality implies  (\ref{ine010s}).

 If the condition (\ref{label1.6}) is satisfied,  by the  point 5 of Lemma \ref{McD}, then we have
    $V_n^2 =  O(n^{1+2\alpha})$  as $n\rightarrow \infty.$
Thus, inequality (\ref{riosbound}) implies that
for any $x>0,$
\begin{eqnarray*}
  \mathbb{P}\left(\| S_n \|_p \geq n x   \right)
    \leq I_1(nx d^{-1/p} )  +  2d\exp \left \{ -c_{p,d} \,\frac{ (x n  )^2}{ n^{1+2\alpha} }
 \right\}.
\end{eqnarray*}
 From the last inequality, we get (\ref{ine011s}).

 If the condition (\ref{label1.7}) is satisfied,  by the point 6 of Lemma \ref{McD}, then we have (\ref{fsdf36}).
Thus, inequality (\ref{riosbound}) implies  that for any $x>0,$
\begin{eqnarray}\label{dfds}
   \mathbb{P}\left(\| S_n \|_p \geq n x   \right)
    \leq   \left\{ \begin{array}{ll}
\displaystyle I_1(nx d^{-1/p} )  +  2 d  \exp\left \{ -c_{p,d} \, \frac{ (x n  )^2}{ n^{1-2\alpha} }\,,
 \right\} \,, & \textrm{if $0\leq \alpha < 1/2,$ }\\
  \\
\displaystyle I_1(nx d^{-1/p} )  +   2 d  \exp\left \{ -c_{p,d} \,\frac{  (x n  )^2}{ \ln n  }
 \right\} \,,  & \textrm{if $ \alpha = 1/2,$ }\\
  \\
\displaystyle I_1(nx d^{-1/p} )  +  2 d  \exp\Big\{ -c_{p,d}(x n   )^2
 \Big\}\,,   & \textrm{if $ 1/2 < \alpha \leq1$.}
\end{array} \right.
\end{eqnarray}  From \eref{dfds}, we obtain the point [iii] of the property.
  \hfill\qed

\subsection{Hoeffding type bound}
The next proposition is an application of Corollary 2.3 of Fan {\it et al.} \cite{FGL12}, which is an extension of Hoeffding's inequality
for super-martingales.
\begin{prop}\label{po71}
Assume that  there exists a positive constant  $A_1$   such that
$$
{\mathbb E} \big[  \big(   G_\varepsilon(\varepsilon)\big)^2\big] \leq
    A_1\, .
 $$
Put
\begin{equation}\label{hoedffineq}
V_n^2 =  2 \sum_{k=2}^n K_{ k,n}^2  (   \tau_{ k} ^2 A_1 +  \xi_k ^2 )  \ \ \ \ and \ \ \ \ \ \delta_n=\max\{  K_{ k,n} \tau_{ k}, K_{ k,n} \xi_k, k=2,\ldots,n\}.
\end{equation}
Then  for any $x,y > 0$,
\begin{equation}\label{Hineq}
{\mathbb P}\Big(\| S_n \|_p \geq  d^{1/p} x\Big)  \leq I_1(x)+   2d \, H_n \left(\frac{x }{2(y+1) \delta_n} , \frac{ V_n  }{(y+1) \delta_n}\right)
   +  2d \,{\mathbb P}\left(\max_{2 \leq k \leq n} G_\varepsilon(\varepsilon_k)  > y \right)   ,
\end{equation}
where
\begin{eqnarray}\label{fkmns}
H_{n}(x,v)=\left\{\left( \frac{v^2}{
x+v^2}\right)^{ x+v^2 }\left( \frac{n}{n-x}\right)^{ n-x }
\right\}^{\frac{n}{n+v^2} }\mathbf{1}_{\{x \leq n\}}\  ,
\end{eqnarray}
with the convention that $(+\infty)^0=1$   \emph{(}which  applies when $x=n$\emph{)}.
In particular, if $$    G_\varepsilon(\varepsilon)   \leq T\ \ a.s.,$$ for a positive constant $T$,
then \eref{Hineq} implies that for any $x  > 0$,
\begin{eqnarray} \label{H}
{\mathbb P}\Big(\| S_n \|_p \geq d^{1/p} x \Big)  \leq I_1(x)+    2d \, H_n \left(\frac{x }{2(T+1) \delta_n} , \frac{ V_n  }{ (T+1) \delta_n}\right) .
\end{eqnarray}
\end{prop}

\noindent\emph{Proof.} We adapt the Corollary 2.3 of Fan {\it et al.}\ \cite{FGL12} with the truncature level  $(y+1) \delta_n$.
By (\ref{ines532}),  we obtain $|M_k^{(i)}| \leq \delta_n (G_\varepsilon(\varepsilon_k)+1) $ for $k\in  [2, n]$ and $i \in [1, d]$. Hence, for any $k\in [2, n]$,
\begin{eqnarray}
 {\mathbb E}\big[ (M_k^{(i)})^2 {\bf 1}_{\{M_k^{(i)} \leq y \delta_n\}}|{\mathcal F}_{k-1}\big]
  &\leq&  2(K_{ k,n}\tau_{ k})^2
 {\mathbb E} \big[  \big( G_\varepsilon(\varepsilon)\big)^2\big] + 2 (K_{ k,n}\xi_k  )^2  \nonumber \\
 & \leq& 2 K_{ k,n}^2  (   \tau_{ k} ^2 A_1 +  \xi_k ^2 ) . \nonumber
\end{eqnarray}
By Corollary 2.3 of Fan  {\it et al.} \cite{FGL12}, it follows   that
\begin{eqnarray}
{\mathbb P}(S_n^{(i)} \geq x/2  )  &\leq&    H_n \left(\frac{x }{2(y+1) \delta_n} , \frac{ V_n   }{(y+1) \delta_n }\right)   + {\mathbb P}\bigg( \max_{2 \leq k \leq n} M_k^{(i)} \geq (y+1) \delta_n \bigg)\, \nonumber  \\
&\leq& H_n \left(\frac{x }{2(y+1) \delta_n} , \frac{ V_n  }{(y+1) \delta_n}\right)
   + {\mathbb P}\left(   \max_{2 \leq k \leq n} G_\varepsilon(\varepsilon_k)  > y \right). \label{fsdv}
\end{eqnarray}
Moreover, the same bound holds for ${\mathbb P}(-S_n^{(i)} \geq x/2) $.
Applying (\ref{fsdv}) to (\ref{f29}), we obtain  the desired inequality. \hfill \qed

\begin{rem}
Using the Remark 2.1 of Fan  {\it et al.}\ \cite{FGL12}, we have  for any $x, v>0$,
\begin{eqnarray}
H_{n}(x, v)&\leq& B(x,v):=\left(\frac{v^2}{x+v^2} \right)^{x+v^2}e^x  \label{freedma1} \\
&\leq& B_1(x,v):= \exp\left\{-\frac{x^2}{2(v^2+\frac{1}{3}x )}\right\} . \label{Bernstein}
\end{eqnarray}
Note that $B(x,v)$ and $B_1(x,v)$ are respectively known as  Bennett's
and Bernstein's bounds. Then, inequality (\ref{H}) also implies the following   Bennett's and Bernstein's bounds. For any $x>0,$ we have
\begin{eqnarray*}
{\mathbb P}\Big(\| S_n \|_p \geq  d^{1/p} x\Big)    &\leq& I_1(x)+    2d \,  B\!\left(\frac{x }{2(T+1) \delta_n} , \frac{ V_n  }{(T+1) \delta_n}\right)  \\&\leq& I_1(x)+    2d \, B_1\!\left(\frac{x }{2(T+1) \delta_n} , \frac{ V_n  }{(T+1) \delta_n}\right).
\end{eqnarray*}
\end{rem}

\section{Moment inequalities}\label{secM}
\setcounter{equation}{0}

\subsection{Marcinkiewicz-Zygmund  type bound}\label{MZsec}
If the martingale differences $(M_i)_{i\geq1}$ have finite $q$th moments ($q\geq 2$), then we have the following Marcinkiewicz-Zygmund type inequality (cf.\ Rio \cite{R09}).
\begin{prop}
 Let $q\geq 2$. Assume that   there exist two positive constants
$B_1(q)$ and $B_2(q)$ such that
\begin{eqnarray}\label{uines41}
&&{\mathbb E}\big[\big(G_{X_1}(X_1)\big)^q \big] \leq B_1(q)\,
\quad   \text{and}  \quad \ \, \ \,
 {\mathbb E}\big[\big(G_\varepsilon(\varepsilon)\big)^q \big] \leq B_2(q)
  \,  .
\end{eqnarray}
Denote
    $$T_n(q)= K_{ 1,n}^2 \big(B_1(q) \big)^{2/q} + (q-1) 2^{2- 2/q}  \sum_{k=2}^n  K_{k, n}^2 \big( \tau_k^q B_2(q)   +  \xi_k^q  \big )^{2/q} .$$
    Then
\begin{equation}\label{ineq4.2}
 {\mathbb E} \| S_n \|_q  \leq  d^{1/q} \sqrt{T_n(q)}.
\end{equation}
\end{prop}

\noindent\emph{Proof.}
  Using Theorem 2.1 of Rio \cite{R09}, we have for any $q \geq 2,$
\begin{eqnarray*}
({\mathbb E} | S_n^{(i)}  |^q )^{2/q} &\leq&  ({\mathbb E} | M_1^{(i)}  |^q )^{2/q} +  (q-1) \sum_{k=2}^n  ({\mathbb E} | M_k^{(i)}  |^q )^{2/q}.
\end{eqnarray*}
Again by (\ref{ines4562}), \eref{ines532}  and (\ref{uines41}),  we deduce that
\begin{eqnarray*}
({\mathbb E} | S_n^{(i)}  |^q )^{2/q} &\leq& \Big(\mathbb{E}[( K_{ 1,n} G_{X_1}(X_1))^q ] \Big)^{2/q} +(q-1) \sum_{k=2}^n\Big(\mathbb{E}[ ( K_{k, n} [\tau_k  G_\varepsilon(\varepsilon_k) + \xi_k ] )^q ] \Big)^{2/q}
\\
&\leq& \Big(\mathbb{E}[( K_{ 1,n} G_{X_1}(X_1))^q ] \Big)^{2/q} +(q-1) \sum_{k=2}^n\Big( 2^{q-1} ( K_{k, n})^q \mathbb{E}[ ( \tau_k  G_\varepsilon(\varepsilon_k))^q  + \xi_k^q    ] \Big)^{2/q}
\\
&\leq&  K_{ 1,n}^2 \big(B_1(q) \big)^{2/q} + (q-1) 2^{2- 2/q}  \sum_{k=2}^n  K_{k, n}^2 \big( \tau_k^q B_2(q)  +  \xi_k^q  \big )^{2/q}\\
&=&   T_n(q) .
\end{eqnarray*}
Using \textcolor{blue}{Jensen's} inequality and the last inequality, we get
\begin{eqnarray*}
 {\mathbb E} \| S_n \|_q  &=& {\mathbb E} \Big(\sum_{i=1}^{d} | S_n^{(i)}  |^q \Big)^{1/q} \leq \Big({\mathbb E}\sum_{i=1}^{d} | S_n^{(i)}  |^q \Big)^{1/q} =   \Big(\sum_{i=1}^{d}{\mathbb E} | S_n^{(i)}  |^q \Big)^{1/q}    \\
 &\leq&  %d^{1/q} \Big[\big(B_1(q) \big)^{2/q}(K_{ 1,n} )^2 +  (q-1) ( B_2(q) )^{2/q}   \sum_{k=2}^n   \big(K_{k, n} \tau_k \big)^2\Big]^{1/2} =
 d^{1/q} \sqrt{T_n(q)}\,,
\end{eqnarray*}
which is the desired inequality.   \hfill\qed

\subsection{von Bahr-Esseen type  bound}\label{VEsec}
When the dominating random variables $G_{X_1}(X_1)$ and $G_\varepsilon(\varepsilon )$ have  a moment of order $q \in [1,2]$,
we have the following von Bahr-Esseen type bound.
\begin{prop}
Let $q \in [1, 2]$. Assume   that
\begin{equation}\label{fnf58}
{\mathbb E} \big[  \big(  G_{X_1}(X_1)\big)^q\big] \leq A_1(q)
\quad  \text{and} \quad
{\mathbb E} \big[  \big(   G_\varepsilon(\varepsilon)\big)^q\big] \leq A_2(q) \, .
\end{equation}
Let
$$V_n(q)= K_{1, n}^q   A_1(q)
         + 2 \sum_{k=2}^n  K_{k, n} ^q    \big( \tau_k^q A_2(q)  + \xi_k^q \big) .$$
Then
\begin{equation}\label{vBE}
 {\mathbb E} \| S_n \|_q  \leq    \big( d \, V_n(q) \big)^{1/q}.
\end{equation}
\end{prop}

\noindent\emph{Proof.} By an argument similar to the proof of  (\ref{fsdsg45}), we have
\begin{eqnarray}
{\mathbb E}  |S_{2,n}^{(i)} |^q   &\leq&       {\mathbb E} \big[  \big(  K_{1, n}    G_{X_1}(X_1)  \big)^q\big]    + 2^{2-q}
    \sum_{k=2}^n  K_{k, n} ^q {\mathbb E} \big[  \big( \tau_k  G_\varepsilon(\varepsilon) + \xi_k \big)^q\big]\,  \nonumber  \\
    &\leq&   K_{1, n}^q  {\mathbb E} \big[  \big(  G_{X_1}(X_1)  \big)^q\big]
         + 2^{q-1} 2^{2-q}
    \sum_{k=2}^n  K_{k, n} ^q {\mathbb E} \big[  \big( \tau_k  G_\varepsilon(\varepsilon))^q + \xi_k  ^q\big]\,   \nonumber  \\
    &\leq&   K_{1, n}^q   A_1(q)
         + 2 \sum_{k=2}^n  K_{k, n} ^q    \big( \tau_k^q A_2(q)  + \xi_k^q \big)    \nonumber  \\
    &  =&  V_n(q) .    \nonumber
\end{eqnarray}
It is easy to see that
\begin{eqnarray*}
 {\mathbb E} \| S_n \|_q  = {\mathbb E} \Big(\sum_{i=1}^{d} | S_n^{(i)}  |^q \Big)^{1/q} \leq \Big({\mathbb E}\sum_{i=1}^{d} | S_n^{(i)}  |^q \Big)^{1/q} =   \Big(\sum_{i=1}^{d}{\mathbb E} | S_n^{(i)}  |^q \Big)^{1/q}   \leq   ( d \, V_n(q) )^{1/q},
\end{eqnarray*}
which gives (\ref{vBE}).     \hfill\qed

\section{Applications}\label{appl}
\setcounter{equation}{0}
\subsection{Stochastic approximation by averaging}\label{SA}
Let us return to the general linear problem  of   Example 3 in Subsection \ref{secexample}. In   Subsection~\ref{SA}, we fix some $p\in[1,\infty]$ and   $\mathcal{X}={\mathcal Y} =\mathbb{R}^d $ equipped with $d(x,x')=\delta(x,x')=\|x-x'\|_p$.
For  linear problem, the central limit theorems for $\overline{X}_n-x^*$   have been well studied by Polyak and Juditsky \cite{PJ92}. In this subsection, we
focus on deviation inequalities for  $\overline{X}_n-x^*$.\medskip

 We first consider the case of Example  3 where $A_i\equiv A$ is deterministic.
Recall that $\alpha \in [0, 1)$, $$Ax^*=B\ \ \ \ \ \ \ \textrm{and} \ \ \ \ \ \ \   X_n  =    X_{n-1} - \frac{\gamma}{n^\alpha} (  A X_{n-1}+\underbrace{\eta_n-B_{n-1}}_{=\varepsilon_n}).$$
From the last line, we deduce that
\begin{eqnarray*}
\mathbb{E} X_n -x^* &=& \mathbb{E} X_{n-1} -  \frac{\gamma}{n^\alpha}A \ \mathbb{E}X_{n-1}  + \frac{\gamma}{n^\alpha}  B -  x^*  \\
&=& \mathbb{E} X_{n-1} -  \frac{\gamma}{n^\alpha} A \ \mathbb{E}X_{n-1}  + \frac{\gamma}{n^\alpha} Ax^*  -  x^* \\
&=& (I_d -\frac{\gamma}{n^\alpha} A )(\mathbb{E}X_{n-1} -x^*) \\
&=&  \prod_{k=2}^{n}(I_d -\frac{\gamma}{k^\alpha} A )(\mathbb{E}X_{1} -x^*) .
\end{eqnarray*}

Thus, we have
\begin{eqnarray*}
\|\mathbb{E} X_n -x^*\|_p  \leq   \|\mathbb{E}X_{1} -x^*\|_p \prod_{k=2}^{n}\rho_k ,
\end{eqnarray*}
where $\rho_n = 1- \frac{\gamma  \lambda_{\min}^{(p)} (A) }{n^\alpha}$.  Using inequality (\ref{dfsfds}), we have $\prod_{k=2}^{n}\rho_k \leq \exp\Big\{ - \frac{\gamma \lambda_{\min}^{(p)} (A)}{ (1-\alpha)^2 } \frac{n-2} {n^\alpha } \Big\},$ which leads to
\begin{eqnarray}
\|\mathbb{E} X_n -x^*\|_p  \ \leq \  \|\mathbb{E}X_{1} -x^*\|_p \exp\Big\{ - \frac{\gamma \lambda_{\min}^{(p)} (A)}{ (1-\alpha)^2 } \frac{n-2} {n^\alpha } \Big\}
\end{eqnarray}
and
\begin{eqnarray*}
\|\mathbb{E} \overline{X}_n -x^*\|_p    \leq  \frac{1}{n} \sum_{k=1}^n\|\mathbb{E} X_k -x^*\|_p    \ \leq \   \frac{C_0}{n} ,
\end{eqnarray*}
where
\begin{eqnarray}\label{fdfs02}
C_0=  \|\mathbb{E}X_{1} -x^*\|_p \bigg( 1+ \sum_{k=2}^{\infty}\exp\Big\{ - \frac{\gamma \lambda_{\min}^{(p)} (A)}{ (1-\alpha)^2 } \frac{k-2} {k^\alpha } \Big\} \bigg).
\end{eqnarray}
Taking  $ f(X_1, X_2, \ldots, X_n)= n \, \overline{X}_n, $  we can see that (\ref{codiMD}).
Clearly, it holds
$$
 \overline{X}_n -x^*      =    \overline{X}_n -\mathbb{E} \overline{X}_n  + \mathbb{E} \overline{X}_n  -x^*,
$$
which implies that
$$
\| \overline{X}_n -x^*  \|_p    \leq  \|  \overline{X}_n -\mathbb{E} \overline{X}_n \|_p  + \| \mathbb{E} \overline{X}_n  -x^*\|_p\leq \|  \overline{X}_n -\mathbb{E} \overline{X}_n \|_p  +   \frac{C_0}{n}.
$$
Hence, we have
\begin{eqnarray*}
{\mathbb P}\Big(\big\|  \overline{X}_n -x^*\big\|_p \geq x\Big) &\leq&  {\mathbb P}\Big(\big\| \overline{X}_n -\mathbb{E} \overline{X}_n\|_p \geq x-\frac{C_0}{n}\Big) \\
&=& {\mathbb P}\Big(\big\| f(X_1, \ldots , X_n) -{\mathbb E}[f(X_1, \ldots , X_n)] \big\|_p \geq nx-C_0 \Big) \\
&=& {\mathbb P}\Big(\big\| S_n  \big\|_p \geq nx-C_0 \Big) .
\end{eqnarray*}
Notice that the condition \eref{label1.5} is satisfied in Example  3.  Thus, the following qualitative inequalities are consequences of  our deviation inequalities.
\begin{itemize}
\item  If \eref{Bernsteinmoment} and \eref{ine35} hold, then
there exist some positive constants $c_{1,p,d}$ and $c_{2,p,d}$,   such that
\begin{equation} \label{fsf52}
 \limsup_{n\rightarrow \infty} \frac{1}{n} \ln {\mathbb P}\Big(\big\|  \overline{X}_n -x^*\big\|_p \geq x\Big) \leq \begin{cases}
  -  c_{1,p,d} \, x       \quad \quad  \text{if\ \ $x \in  (1, \infty)$}\\
  \vspace{-0.2cm}\\
    -  c_{2,p,d} x^2     \quad  \quad  \text{if\ \ $x \in (0, 1]$.}
  \end{cases}
\end{equation}
This follows from the point  [i] in Proposition \ref{prof3.1}.
\item If \eref{laplace2}  and \eref{ifnelh5} hold  for some $q \in (0,1)$, then
there exists  a positive constant $c_{p,d}$ such that for any $x >0 ,$
\begin{eqnarray}\label{fsf53}
 \limsup_{n\rightarrow \infty} \frac{1}{n^q } \ln {\mathbb P}\Big(\big\|  \overline{X}_n -x^*\big\|_p \geq x\Big) \leq  - c_{p,d}\, x^q   .
\end{eqnarray}
This follows from the point  [i] in Proposition \ref{findsa}.
\item  If \eref{ifgnelh5}  and \eref{fnafdsa58} hold  for some $q\geq1$,  then there exists a positive
constant $c_{p,d}$ such that for any $x >0,$
\begin{eqnarray}\label{fsf54}
\limsup_{n\rightarrow \infty} n^{q-1} {\mathbb P}\Big(\big\|  \overline{X}_n -x^*\big\|_p \geq x\Big) \leq \frac{c_{p,d}}{ x^q} \, .
\end{eqnarray}
 This follows from the points [i] in Proposition \ref{mineq} (case $q \geq 2$) and  Proposition \ref{VBEI} (case $q \in [1, 2))$.
\end{itemize}
And for the moment bounds of $S_n$:
\begin{itemize}
\item  If \eref{uines41} holds  for some $q\geq 2$, then,   by \eref{ineq4.2} and the point 4 of Lemma \ref{McD},
\begin{equation}\label{ROSEN}
{\mathbb E}\| \overline{X}_n -x^*  \|_q    \leq  \frac1n {\mathbb E}\| S_n \|_q  +   \frac{C_0}{n} \leq  \frac{c_{p,d}}{\sqrt{n}} .
\end{equation}
\item  If \eref{fnf58} holds for some $q \in [1,2]$, then,  by \eref{vBE} and $ V_n(q)= O(n)$,
\begin{equation} \label{fsf55}
{\mathbb E}\| \overline{X}_n -x^*  \|_q \leq {\mathbb E}\|  \overline{X}_n -\mathbb{E} \overline{X}_n \|_q  +   \frac{C_0}{n}  \leq \frac1n {\mathbb E}\| S_n \|_q  +   \frac{C_0}{n}\leq \frac{c_{p,d}}{ n^{1-1/q}} .
\end{equation}
\end{itemize}
\begin{rem}
Let us make some comments on the performances of uniform averaging $\overline{X}_n$, final iterate $X_n$ and suffix averaging $ \hat{X}_n= \frac{2}{n}\sum_{i=[n/2]}^{n}X_i$.
\begin{enumerate}
  \item Clearly, if $\overline{X}_n -x^*$ is replaced by $ X_n-x^*$, then   inequalities \eref{fsf52}-\eref{ROSEN} hold true,  with $C_0$
  replaced by
  $$C_n=   n \|\mathbb{E}X_{1} -x^*\|_p \exp\Big\{ - \frac{ \lambda_{\min}^{(p)} (A)}{ (1-\alpha)^2 } \frac{n-2} {n^\alpha } \Big\}  $$
which is smaller than $C_0$ defined by \eref{fdfs02} for any $n$ large enough.  Hoverer, this does not improve the convergence rates for the bounds \eref{fsf52}-\eref{ROSEN}. Thus uniform averaging $\overline{X}_n$ and  final iterate $X_n$  have almost the same performance for estimating $x^*$ in a long time view.
  \item When uniform averaging $\overline{X}_n$ is replaced by suffix averaging $\hat{X}_n$, the inequalities \eref{fsf52}-\eref{ROSEN} remain valid. %, instead of uniform averaging $\overline{x}_n$.
\end{enumerate}
\end{rem}

%\begin{rem}
%Let $x_1'$ and $x_1''$ be two different point or random  points.
%If $\overline{x}_n -x^*$ is replaced by $ x'_n-x''_n$, then the inequalities \eref{fsf52}-\eref{ROSEN} hold true.
%\end{rem}
 
Let us now focus on the case where $A_i$ is stochastic. Recall that $\alpha \in [0, 1)$,
$$ x^* =\arg\min_{x\in\mathcal{C}}\|Ax-B\|_2   \ \ \ \ \ \ \ \textrm{and} \ \ \ \ \ \ \   X_n  =  \Pi_\mathcal{C} \left[  X_{n-1} - \frac{\gamma}{n^\alpha} (  A_{n-1} X_{n-1}+\eta_n-B_{n-1}) \right].$$
For the sake of simplicity, assume that $X_1\in\mathcal{C}$,  and  that the condition (\ref{Bernsteinmoment}) is satisfied. Under our assumptions, with moreover $\alpha\in[1/2,1)$, and assuming that $\eta_n$ and $B_{n-1}$ have moments of order $4$, Theorem 3 of~\cite{MB2011} leads to
$$ \mathbb{E}  \|\overline{X}_n -x^*\|_2   = \frac{C_0}{\sqrt{n}} $$
for some $C_0>0$.
Combining this with the point [i] in Proposition~\ref{prof3.1},  we obtain
   for any $x>0,$
\begin{eqnarray}
  \limsup_{n \rightarrow \infty} \frac{1}{n}\ln \mathbb{P}\left(\| \overline{X}_n -x^* \|_2 \geq   x   \right)
    \ \leq \    - \,  c_{2,d}  \, \big( x \mathbf{1}_{\{x\geq 1\}} + x^2 \mathbf{1}_{\{0< x  <  1\}} \big)  ,
\end{eqnarray}
  from which we get   for any $\delta \in (0, 1),$
$$
\mathbb{P}\left( \|\overline{X}_n -x^*\|_2 \  \leq  \
\frac{C_0 +\sqrt{\frac{ 2 }{c_{2, d}}\ln\frac{1}{\delta}} } { \sqrt{n}}
  \ \right) \geq 1-\delta
$$
 for all $n$ large enough.

\subsection{Empirical risk minimization}\label{ERM}

It has been shown in the past how Bernstein type inequalities allow to control the error of the empirical risk minimizer, and thus to perform model selection for time series~\cite{M00,AW11,ALW13,HS14,ADF19,ABDG20}. These results are available under restrictive assumptions. For example,~\cite{AW11,ALW13,ABDG20} focus on stationary series. In \cite{ADF19}, nonstationary Markov chains as in~\eqref{Mchain} are considered, under the restriction that $\rho_n \leq \rho < 1$ and $\tau_n \leq \eta$ in the conditions~\eqref{contract} and~\eqref{c2}. Our new Bernstein type bound, Proposition~\ref{prof3.1}, allows to extend these results to a more general setting.

The context is as follows. For simplicity, here, $\mathcal{X}$ will be a Banach space with norm $\|\cdot \|_{\mathcal{X}}$. Assume we have a parameter set $\Theta$ and a family of functions $f_n(\theta,x)$ of $\theta\in\Theta$ and $x\in \mathcal{X}$, with $f_n(\theta,0)=0$. We observe
that $X_1,\dots,X_n$ satisfy~\eqref{Mchain}, with $F_n(x,y)=f_n(\theta^0,x)+y$ for some unknown $\theta^0\in\Theta$. Of course, the distributions  of $X_1$ and of the $\varepsilon_n$'s  are  also unknown.

\begin{rem}
 We review here some examples studied in the aforementioned references. In the case $\mathcal{X} = \mathbb{R}^d$,~\cite{ABDG20} studied functions of the form
 $$ f_n(\theta,x) = \theta x $$
 where $\theta$ is some $d\times d$ matrix. Note that in this case, the model does actually not depend on $n$. On the other hand,~\cite{ADF19} considered a $T$-periodic version of these functions: for $\theta = (A_1|\dots|A_T)$ where each $A_t$ is a $d\times d$ matrix, they used
 $$ f_n(\theta,x) = A_{n ({\rm mod}\, T)} x. $$
 Other examples include nonlinear autoregression with neural networks~\cite{AW11}.
\end{rem}

Let $\ell:\mathcal{X} \rightarrow [0,+\infty) $ be a function with $\ell(0)=0$, it is usually refered to as the loss function. We will measure the performance of a predictor through its risk:
$$ R_n(\theta) = \frac{1}{n-1} \sum_{k=2}^n \mathbb{E}\big[\ell\big(X_k - f_k(\theta,X_{k-1}) \big) \bigr]. $$
A classical loss function is simply given by the norm $\ell(x) = \|x\|_{\mathcal{X}} $ but other examples can be used, for example~\cite{ALW13} used quantile losses in the case $\mathcal{X}=\mathbb{R}$. Our objective will be to estimate the minimizer $\theta^*$ of $R_n$. Under suitable assumptions, $\theta^*=\theta^0$: this is for example the case when $\mathcal{X}$ is actually a Hilbert space, $\ell(x) = \|x\|_{\mathcal{X}}^2$ and the $\varepsilon_n$ are centered with $\mathbb{E} \|\varepsilon_n\|_{\mathcal{X}}^2 <\infty$. However, this has no reason to be true in general, and it is important to note that if the objective is to minimize the loss of the predictions, to estimate $\theta^*$ is more important than to estimate $\theta^0$.  We define the ERM estimator $\hat{\theta}$ (for Empirical Risk Minimizer) by
$$
\hat{\theta}=\arg\min_{\theta\in\Theta} r_n(\theta),
\text{\ \ \ \  where } r_n(\theta) = \frac{1}{n-1} \sum_{k=2}^n \ell\big(X_k - f_k(\theta,X_{k-1}) \big).
$$
\begin{Defi}
 \label{defi.covering}
 Define the covering number $\mathcal{N}(\Theta,\epsilon)$ as the cardinality of the smallest set $\Theta_{\epsilon} \subset \Theta$ such that for any  $ \theta\in\Theta$, there exists a $ \theta_{\epsilon} \in\Theta_{\epsilon}$ such that
 $$ \sup_{k\in\{2,\dots,n\}} \sup_{x\in\mathcal{X}}
 \frac{\|f_k(\theta,x)-f_k(\theta_\epsilon, x)\|_{\mathcal{X}} }{ \| x\|_{\mathcal{X}} }
  \leq \epsilon .$$
 Define the entropy of $\Theta$ by  $\mathcal{H} (\Theta,\epsilon) = 1\vee \ln \mathcal{N}(\Theta,\epsilon)$.
\end{Defi}

Examples of computation of $\mathcal{H} (\Theta,\epsilon)$ for some models can be found in references~\cite{ADF19,ABDG20}. In most classical examples, $\mathcal{H} (\Theta,\epsilon)$ is roughly in $1\vee [D \ln (1+C/\epsilon)],$ where $D$ is the dimension of $\Theta$ and $C>0$ is some constant.
\begin{prop}
Let $X_1,\dots,X_n$ satisfy~\eqref{Mchain},~\eqref{contract} and~\eqref{c2} with $(\tau_n) $ and $(\rho_n)$ satisfying \eqref{label1.6} with $\alpha \in (0, 1/2)$  and  $d(x , y)=\delta(x , y )=\| x -y  \|_{\mathcal{X}}$ (note that $   \xi_n \equiv 0$ in this case). Assume that~\eqref{Bernsteinmoment} is satisfied, and that $\ell$ is $L$-Lipschitz, that is for any $(x,y)\in\mathcal{X}^2$,
$$| \ell (x) - \ell(y) | \leq L    \|x -y \|_{\mathcal{X}}    . $$
 Assume also that all the functions in the model are $\lambda$-Lipschitz: for any $(x,y)\in\mathcal{X}^2$,   any $k\in\mathbb{N}$ and  any $\theta\in\Theta$, $$   \| f_k(\theta,x) - f_k(\theta,y) \|_{\mathcal{X}  }\leq \lambda \|x-y\|_{\mathcal{X}} ,$$
and that $\mathcal{H}(\Theta,1/(Ln))\leq D \ln  n  $ for some constant $D$.
For $n$ large enough, we have  for any $\eta \in (0, 1)$,
\begin{equation*}
R_n(\hat{\theta})
\leq
 \min_{\theta\in\Theta } R_n(\theta) + C_1 \sqrt{ \frac{D\ln n}{n^{1-2\alpha}}} + C_2 \frac{1+\ln\big(\frac{1}{\eta}\big)}{\sqrt{n^{1-2\alpha}}}
\end{equation*}
with probability at  least $1-\eta,$
where $C_1$ and $C_2$ are constants that depend only on $\lambda$, $L$, and the constant $c_{p,d}$ in the proof of Proposition~\ref{prof3.1}.
\end{prop}

In particular, for $\alpha=0$, we recover bounds that are similar to the ones in~\cite{ADF19,ABDG20}.

In the case where several models are available and one doesn't know which one contains the truth, the previous result can be used to perform model selection. We refer the reader to~\cite{ADF19} for example for details on this classical construction.

\noindent\emph{Proof.}
Fix $\theta\in\Theta$ and consider the random variable $S_n = g(X_1,\dots,X_n)-\mathbb{E}[g(X_1,\dots,X_n)]$, where
$$ g(x_1,\dots,x_n) = \frac{1}{L(\lambda+1)} \sum_{k=2}^n \ell\left(x_k - f_k(\theta,x_{k-1}) \right). $$
Note that
\begin{align*}
| g(x_1,\dots,x_n) - g(x_1,\dots,x_k',\dots,x_n) |
    & \leq \frac{\left| \ell \left(x_{k+1} - f_{k+1}(\theta,x_{k}) \right) - \ell \left(x_{k+1} - f_{k+1}(\theta,x_{k}') \right)\right|}{L(\lambda+1)}
    \\
      & \quad \quad  + \frac{ \left|  \ell \left(x_k - f_k(\theta,x_{k-1}) \right) - \ell \left(x_k' - f_k(\theta,x_{k-1}) \right) \right|}{L(\lambda +1)}
        \\
   & \leq  \frac{\|f_{k+1}(\theta,x_{k})-f_{k+1}(\theta,x_{k}')\|_{\mathcal{X}}  + \|x_k-x_k'\|_{\mathcal{X}} }{\lambda+1}     \\
   & \leq       \|x_k - x_k'\|_{\mathcal{X}}   ,
\end{align*}
which means that $g$ is separately Lipschitz. So  we apply~\eqref{fsdfs} in the proof of Proposition~\ref{prof3.1}, that is, for any $t\in[0,\delta_n^{-1})$,
$$ \mathbb{E}[\exp\left\{ \pm t S_n \right\}] \leq \exp\left\{\frac{t^2 V_{n}}{2-2t\delta_n }\right\} .$$
Note that $S_n = \frac{n-1}{L(1+\lambda)}\left( r_n(\theta)- \mathbb{E}[r_n(\theta)]\right)$, and $ R_n(\theta) = \mathbb{E}[r_n(\theta)] $.
Set $s= {t(n-1)}/{L(1+\lambda)}$. We obtain that, for any $s\in[0,{\delta_n^{-1}(n-1)}/{L(1+\lambda)})$,
\begin{equation}\label{proof.ERM.moment.expo}
\mathbb{E}[\exp\left\{ \pm s ( r_n(\theta)- R_n(\theta)) \right\}]
\leq \exp\left\{ \frac{s^2(1+\lambda)^2 L^2 \frac{V_{ n }}{n-1}}{2(n-1)-2s(1+\lambda)\delta_{n} L} \right\}.
\end{equation}

Fix now $\epsilon>0$ and a set $\Theta_\epsilon\subset\Theta$ as in Definition~\ref{defi.covering}. For any $\theta\in\Theta_\epsilon$, we have $\theta\in\Theta$ and so~\eqref{proof.ERM.moment.expo} holds. Then, for any $s\in[0,{\delta_n^{-1}(n-1)}/{L(1+\lambda)})$ and any $x>0$, we have
\begin{align}
 \mathbb{P}\left(\sup_{\theta\in\Theta_\epsilon}|r_n(\theta)-R_n(\theta)| > x \right)
& \leq \sum_{\theta\in\Theta_\epsilon} \mathbb{P}\left(|r_n(\theta)-R_n(\theta)| > x \right)
\nonumber
\\
& \leq \sum_{\theta\in\Theta_\epsilon} \mathbb{E}[\exp\left\{ s |r_n(\theta)-R_n(\theta)| - s x\right\}]
\nonumber
\\
& \leq 2 \mathcal{N}(\Theta,\epsilon) \exp\left\{ \frac{s^2(1+\lambda)^2 L^2 \frac{V_{ n }}{n-1}}{2(n-1)-2s(1+\lambda)\delta_{n} L} -sx \right\}.
\label{equa.dev.1}
\end{align}
Thanks to the definition of $\Theta_\epsilon$, for any $\theta\in\Theta$ there is a $\theta_{\epsilon}$ such that
 $$ \sup_{i\in\{2,\dots,n\}} \sup_{x\in\mathcal{X}}  \frac{\|f_i(\theta,x)-f_i(\theta_\epsilon,x)\|_{\mathcal{X}} }{\|x\|_{\mathcal{X}}  }      \leq \epsilon .$$
So
$$
\left| \ell(X_k-f_k(\theta_\epsilon,X_{k-1})) -  \ell(X_k - f_k(\theta,X_{k-1}) )\right|
\leq
L \|f_k(\theta_\epsilon,X_{k-1})   - f_k(\theta,X_{k-1})\|_{\mathcal{X}}
\leq L \epsilon \|X_{k-1}\|_{\mathcal{X}}
$$
and as a consequence,  we obtain
\begin{equation}\label{fsggss01}
  |r_n(\theta)-r_n(\theta_\epsilon)| \leq   L \epsilon \cdot\frac{\displaystyle 1}{n-1}\sum_{k=1}^{n-1}
   \|X_k\|_{\mathcal{X}}
\end{equation}
and
\begin{equation} \label{fsggss02}
|R_n(\theta)-R_n(\theta_\epsilon)| \leq L \epsilon \cdot\frac{\displaystyle 1}{n-1}\sum_{k=1}^{n-1}\mathbb{E}
 \|X_k\|_{\mathcal{X}}.
\end{equation}
Using Proposition~\ref{prof3.1} with $f(X_1,\dots,X_n)=\sum_{k=1}^{n-1} \|X_k\|_{\mathcal{X}}, $ we get for any $y>0$ and any $u\in[0,\delta_n^{-1})$,
\begin{align}
\mathbb{P}\Big(\sum_{k=1}^{n-1}\|X_k\|_{\mathcal{X}} > \sum_{k=1}^{n-1}\mathbb{E}\|X_k\|_{\mathcal{X}}  + y \Big)
& \leq \mathbb{E}\exp\bigg\{u\Big(\sum_{k=1}^{n-1}\|X_k\|_{\mathcal{X}} - \sum_{k=1}^{n-1}\mathbb{E}\|X_k\|_{\mathcal{X}}  - y\Big) \bigg\}
\nonumber
\\
& \leq \exp\left\{\frac{u^2 V_{ n }}{2\left(1-u\delta_n \right)} - uy \right\}.
\label{equa.dev.2}
\end{align}
From now, let us use the short notation $z_{n} = \sum_{k=1}^{n-1}\mathbb{E} \|X_k\|_{\mathcal{X}} $ and consider the ``favorable'' event
$$
\mathcal{E}= \left\{ \sum_{k=1}^{n-1}\|X_k\|_{\mathcal{X}}  \leq z_{n} + y \right\}
  \bigcap \left\{\sup_{\theta\in\Theta_\epsilon}|r_n(\theta)-R_n(\theta)| \leq x  \right\} .
$$
On $\mathcal{E}$, by \eref{fsggss01} and  \eref{fsggss02}, we have
\begin{align*}
R_n(\hat{\theta})
& \leq R_n(\hat{\theta}_\epsilon) +\epsilon L\frac{z_{n} }{n-1}   \leq r_n(\hat{\theta}_\epsilon) + x  +\epsilon L \frac{z_{n} }{n-1}
\\
& \leq r_n(\hat{\theta}) + x + \epsilon L \left[ 2  \frac{z_{n} }{n-1}+ \frac{y}{n-1}\right]
\\
& = \min_{\theta\in\Theta } r_n(\theta) + x + \epsilon L   \frac{2 z_{n}+y }{n-1}
\\
& \leq \min_{\theta\in\Theta_\varepsilon } r_n(\theta) + x + \epsilon L   \frac{2 z_{n}+y }{n-1}
\\
& \leq \min_{\theta\in\Theta_\epsilon } R_n(\theta) + 2x + \epsilon L   \frac{2 z_{n}+y }{n-1}
\\
& \leq \min_{\theta\in\Theta } R_n(\theta) + 2x + \epsilon L  \frac{3 z_{n}+y }{n-1}  .
\end{align*}
In particular, the choice $\epsilon = 1/(Ln)$ ensures:
\begin{equation}
\label{almost-done}
R_n(\hat{\theta})
\leq
 \min_{\theta\in\Theta } R_n(\theta) + 2x +  \frac{3 z_{n}+y }{n(n-1)}.
\end{equation}

Inequalities~\eqref{equa.dev.1} and \eqref{equa.dev.2} lead to
\begin{equation}
\mathbb{P}\left( \mathcal{E}^c\right)\leq  \exp\left\{\frac{u^2 V_{ n }}{2\left(1-u\delta_n \right)} - uy \right\}
+2 \mathcal{N}(\Theta,\frac{1}{Ln})
 \exp\left\{ \frac{s^2(1+\lambda)^2 L^2 \frac{V_{ n }}{n-1}}{2(n-1)-2s(1+\lambda)\delta_{n} L} -s x \right\}.
\end{equation}
As it was explained in the proof of Proposition~\ref{prof3.1}, we get $\delta_n= O (n^\alpha)$ and $V_n= O (n^{1+2\alpha})$. So, letting $c_{p,d}$ be as in the proof of Proposition~\ref{prof3.1}, we get
\begin{equation}\label{proba-e}
\mathbb{P}\left( \mathcal{E}^c\right)\leq  \exp\left\{\frac{u^2 c_{p,d}n^{1+2\alpha}}{2\left(1-u c_{p,d}  n^{\alpha} \right)} - uy \right\}
+2 \mathcal{N}(\Theta,\frac{1}{Ln}) \exp\left\{ \frac{s^2(1+\lambda)^2 L^2 c_{p,d}  n^{2\alpha} }{2(n-1)-2s(1+\lambda) c_{p,d} n^{\alpha} L} -sx\right\}.
\end{equation}
Fix $\eta \in (0, 1)$ and put
$$ x =  \frac{s(1+\lambda)^2 L^2 c_{p,d} n^{2\alpha} }{2(n-1)-2s(1+\lambda)c_{p,d} n^{\alpha} L} + \frac{\mathcal{H} (\Theta,\frac{1}{Ln}) + \ln\big(\frac{4}{\eta}\big)}{s} $$
and
$ y =  \frac{\ln\big(\frac{2}{\eta}\big)}{u} + \frac{u c_{p,d} n^{1+2\alpha}}{2(1-u c_{p,d} n^{\alpha} )}.$
Note that, plugging $x$ and $y$ into~\eqref{proba-e}, those choices ensure $\mathbb{P}(\mathcal{E}^c)\leq \eta/2+\eta/2=\eta$, while \eqref{almost-done} becomes:
\begin{eqnarray*}
R_n(\hat{\theta})
 &\leq&
 \min_{\theta\in\Theta } R_n(\theta) +  \frac{s(1+\lambda)^2 L^2 c_{p,d} n^{2\alpha} }{(n-1)-2s(1+\lambda)c_{p,d} n^{\alpha} L} + 2 \frac{D \ln(n) + \ln \big(\frac{4}{\eta}\big)}{s}
 \\
&&\    +  \frac{3 z_{n}}{n(n-1)}
 + \frac{\ln\big(\frac{2}{\eta}\big)}{u n (n-1)} +  \frac{u c_{p,d} n^{1+2\alpha}}{2(1-u c_{p,d} n^{\alpha} ) n (n-1)}.
\end{eqnarray*}
The final steps are to choose $u\in[0,\delta_n^{-1})$, $s\in[0,{\delta_n^{-1}(n-1)}/{L(1+\lambda)})$ and to provide an upper bound on $z_n$. First, put $u=1/(2 c_{p,d} n^\alpha )$ and
$$s= \frac{n-1}{(1+\lambda) L} \sqrt{\frac{2 D\ln(n)}{n^{1+2\alpha} }}. $$
Note that we always have $u < \delta_n^{-1}$. Moreover, we have $s n^{\alpha} = o(n) $, so for $n$ large enough, the condition on $s$ is satisfied too. Thus, for $n$ large enough, there are constants $C_1$ and $C_2$ such that
\begin{equation}\label{ineq4.15}
R_n(\hat{\theta})
\leq
 \min_{\theta\in\Theta } R_n(\theta) + C_1 \sqrt{ \frac{D\ln n}{n^{1-2\alpha}}} + C_2 \frac{1+\ln\big(\frac{1}{\eta}\big)}{\sqrt{n^{1-2\alpha}}} + \frac{3z_n}{n(n-1)}.
\end{equation}
Let us now present an upper bound of $ z_n = \sum_{k=1}^{n-1}\mathbb{E}  \|X_k\|_{\mathcal{X}}$.
 Recall that $  \mathbb{E} \|X_1\|_{\mathcal{X}}  =  \int \|x-0\|_{\mathcal{X}}  \mathbb{P}_{X_1}(dx)  = G_{X_1}(0)  $ and $    \mathbb{E} \|\varepsilon \|_{\mathcal{X}}  =  \int \|x-0\|_{\mathcal{X}} \mathbb{P} _{\varepsilon }(d x)    = G_{\varepsilon }(0). $  Then, by the point 5 of Lemma \ref{McD} and inequality (\ref{dfsfds}),  we have
 \begin{align*}
 \mathbb{E} \|X_n\|_{\mathcal{X}}  & = \mathbb{E} \|f_{n-1} (\theta, X_{n-1}) + \varepsilon_n \|_{\mathcal{X}}
  \leq \mathbb{E} \|f_{n-1} (\theta, X_{n-1})-f_{n-1} (\theta, 0) \|_{\mathcal{X}}  + \mathbb{E} \|\varepsilon_n \|_{\mathcal{X}}  \\
 & \leq \rho_n \mathbb{E} \|X_{n-1}\|_{\mathcal{X}}  + G_\varepsilon(0) \leq  \cdot\cdot\cdot \leq  \rho_n\cdot\cdot\cdot \rho_2 \mathbb{E} \|X_{1}\|_{\mathcal{X}}  + K_{2,n} \mathbb{E} \|\varepsilon \|_{\mathcal{X}}  \\
 &\leq C_3 n^\alpha G_{\varepsilon}(0)  + \exp\Big\{ - \frac{(n-1)\rho }{ (1-\alpha)^2 n^{\alpha}}  \Big\}  G_{X_1}(0).
 \end{align*}
 Therefore, it holds $$z_n \leq C_4 n^{1+\alpha}. $$
   Applying the last inequality to \eref{ineq4.15}, we get
   \begin{eqnarray*}
R_n(\hat{\theta})
&\leq&  \min_{\theta\in\Theta } R_n(\theta) + C_1 \sqrt{ \frac{D\ln n}{n^{1-2\alpha}}} + C_2 \frac{1+\ln\big(\frac{1}{\eta}\big)}{\sqrt{n^{1-2\alpha}}} + \frac{3C_4 n^{1+\alpha}}{n(n-1)} \\
&\leq& \min_{\theta\in\Theta } R_n(\theta) + C_5 \sqrt{ \frac{D\ln n}{n^{1-2\alpha}}} + C_2 \frac{1+\ln\big(\frac{1}{\eta}\big)}{\sqrt{n^{1-2\alpha}}} ,
\end{eqnarray*}
which gives the desired result.
 \hfill \qed

\section*{Acknowledgements}
This work has been partially supported by the National Natural Science Foundation
of China (Grant nos.\,11601375 and 11971063). This work has also been   funded by CY Initiative of Excellence
(grant ``Investissements d'Avenir" ANR-16-IDEX-0008),
Project ``EcoDep" PSI-AAP2020-0000000013.
%The author deeply indebted to the editor and the anonymous referee  for their helpful comments.
%This work has been partially supported by the National Natural Science Foundation
%of China (Grant nos.\,11601375 and 11971063).

\end{document}